%% file: text.tex
\documentclass[a4paper,reqno,oneside]{article}
\newcommand{\nameofthesis}{Finitary Boolean functions}

\usepackage[activate={true,nocompatibility},final,tracking=true,kerning=true,spacing=true,factor=1100,stretch=10,shrink=10]{microtype}
\microtypecontext{spacing=nonfrench}
\usepackage{amsmath}%
\usepackage{amsfonts}%
\usepackage{amssymb}%
\usepackage{amsthm}
\usepackage{graphicx}
\usepackage{mathtools}
\usepackage[utf8]{inputenc}
\usepackage{color,soul}
\usepackage{varioref}

\usepackage{stmaryrd}

\usepackage{csquotes}
\usepackage[natbib, backend=biber, style=alphabetic, backref=true]{biblatex}
\bibliography{bibliography}

\usepackage{pgf,tikz}
\usepackage{mathrsfs}
\usetikzlibrary{arrows}

\usepackage{epigraph}
\usepackage{listings}
\usepackage{mathrsfs}
\usepackage{dsfont}

\usepackage{float}

\usepackage{verbatim}

\usepackage[a4paper, total={5.5in, 8in}]{geometry}

\usepackage[pdfauthor={Vilhelm Agdur},%
pdftitle={\nameofthesis},%
pdftex]{hyperref}

\newcommand{\smartref}[2]{#1 \ref{#2}}

\newtheorem{theorem}{Theorem}[section]
\newtheorem{lemma}[theorem]{Lemma}

\newtheorem{prop}[theorem]{Proposition}
\newtheorem{cor}[theorem]{Corollary}

\theoremstyle{definition}
\newtheorem{definition}[theorem]{Definition}
\newtheorem{example}[theorem]{Example}
\newtheorem{exercise}[theorem]{Exercise}

\newtheorem{remark}[theorem]{Remark}

\numberwithin{equation}{section}

\usepackage[toc,page]{appendix}

\newcommand{\abs}[1]{{\left\lvert#1\right\rvert}}

\newcommand{\floor}[1]{\left\lfloor#1\right\rfloor}

\newcommand{\R}{\mathbb{R}}

\newcommand{\N}{\mathbb{N}}
\newcommand{\Z}{\mathbb{Z}}

\newcommand{\F}{\mathcal{F}}

\newcommand{\E}[1]{\mathbb{E}\left[ #1 \right]}
\newcommand{\Var}[1]{\operatorname{Var}\left(#1\right)}
\newcommand{\Prob}[1]{\mathbb{P}\left( #1 \right)}

\newcommand{\ind}[1]{\mathds{1}_{\left\{#1\right\}}}
\newcommand{\given}{\ \middle|\ }

\newcommand{\Pset}{\mathcal{P}}

\newcommand{\inpr}[2]{\langle#1,#2\rangle}
\newcommand{\norm}[1]{\left\lVert#1\right\rVert}

\newcommand{\symdiff}{\mathbin{\triangle}}

\title{\nameofthesis}
\author{Vilhelm Agdur}
\date{May 22, 2019\footnote{Originally on this date as a Master's Thesis at Gothenburg University, with Jeffrey Steif as advisor. Added some figures when publishing on the arXiv on June 9 2019.}}

\begin{document}
\lstset{language=matlab}

\definecolor{light-red}{rgb}{0.8,0,0}
\definecolor{light-blue}{rgb}{0,0,0.8}
\definecolor{pale-pink}{rgb}{1,0.95,0.95}
\definecolor{cream-white}{rgb}{1,0.99,0.99}
\definecolor{pure-white}{rgb}{1,1,1}
\definecolor{pale-green}{rgb}{0.6,0.89,0.68}
\maketitle

\begin{abstract}
    We study functions on the infinite-dimensional Hamming cube $\{-1,1\}^\infty$, in particular Boolean functions into $\{-1,1\}$, generalising results on analysis of Boolean functions on $\{-1,1\}^n$ for $n\in\N$. The notion of noise sensitivity, first studied in \cite{BKS_theorem}, is extended to this setting, and basic Fourier formulas are established. We also prove hypercontractivity estimates for these functions, and give a version of the Kahn-Kalai-Linial theorem giving a bound relating the total influence to the maximal influence.
    
    Particular attention is paid to so-called finitary functions, which are functions for which there exists an algorithm that almost surely queries only finitely many bits. Two versions of the Benjamini-Kalai-Schramm theorem characterizing noise sensitivity in terms of the sum of squared influences are given, under additional moment hypotheses on the amount of bits looked at by an algorithm. A version of the Kahn-Kalai-Linial theorem giving that the maximal influence is of order $\frac{\log(n)}{n}$ is also given, replacing $n$ with the expected number of bits looked at by an algorithm.
    
    Finally, we show that the result in \cite{schramm_steif_2010} that revealments going to zero implies noise sensitivity also holds for finitary functions, and apply this to show noise sensitivity of a version of the voter model on sufficiently sparse graphs.
\end{abstract}

\tableofcontents

\input{1_definition_and_elementary_properties.tex}
\input{2_noise_sensitivity_basics.tex}
\input{3_0_finite_approximations_and_hypercontractivity.tex}
\input{4_0_influences_and_BKS.tex}

\input{5_revealment_theorem.tex}
\input{6_example_voter_model.tex}



\printbibliography
\end{document}

%% file: 1_definition_and_elementary_properties.tex
\section{Introduction}

Suppose a democratic country has $n$ people living in it, and they are holding a referendum about whether hot dogs or hamburgers are better. Each voter randomly chooses a preference, with each being equally likely, on the day she turns five years old. If we encode the two options as $+1$ and $-1$, we can mathematically formalise their election as a function from $\{-1,1\}^n$ into $\{-1,1\}$. 

Now, suppose the people who count the votes do not particularly care about such a frivolous referendum, and so they are careless in counting the vote. In particular, with probability $\epsilon$, they record a completely random vote, instead of the vote that the citizen made. The question we are interested in is this: Will this random error change the outcome of the election?

The election function $f_n: \{-1,1\}^n\to \{-1,1\}$ is an example of a Boolean function, and the question we have asked is roughly whether this function is noise sensitive or not, as originally defined in \cite{BKS_theorem}. This is a concept that has wide applications, spreading far from our starting point of culinary democracy.

A staggering array of things can be studied using Boolean functions: elections with various voting rules with Condorcet's paradox and Arrow's theorem\cite{KALAI_2002}, properties of random graphs\cite{friedgut_1999, lubetzky_steif_2015}, percolation and other models from statistical mechanics\cite{garban_steif_2015}, various problems of computational complexity and learning in computer science\cite{kalai_safra_2005, odonnell_wimmer_2009, friedgut_1999}, and so on.

It is also a concept that can be generalised in various directions. Aside from varying the probability of $1$ and $-1$, one can study functions on a generalised product domain $\Omega^n$ for some finite set $\Omega$ other than $\{-1,1\}$.

If one views $\{-1,1\}^n$ as an $n$-dimensional Hamming cube, one can see the error in the votes as being caused by a random walk on this cube. This point of view turns out to be very enlightening, and opens up connections to the study of functions on $\R^n$ with Gaussian measure and other generalisations.

The general idea of these generalisations is to replace the random walk on the cube with a random walk on some other structure or some other random process on the cube, leading to results on noise sensitivity with respect to general Markov chains\cite{forsstrom_2015}, exclusion sensitivity\cite{Broman2013}, functions on Gaussian space, for families of log-concave measures, and for functions on Cayley and Schreier graphs. Interestingly enough, all of the latter three and the Hamming cube case can be treated in one unified way using the semigroups of the random walks.\cite{bouyrie_2018}

So, given this swathe of applications and generalisations, what is it we will do? Consider the following question:

Suppose you are given an infinite sequence of $1$s and $-1$s, and told that these are the steps of a one-dimensional random walk started at zero. You are asked to tell whether this one-dimensional random walk hit $10$ before it hit $-10$.

This is very nearly a Boolean function as studied above, except it depends on infinitely many bits. It is, however, close to being one, since it will almost surely hit $10$ or $-10$ in finitely many steps, and so one does not have to actually see infinitely many bits.

One gets the feeling that it should be possible to use the tools of standard analysis of Boolean functions to answer questions about this function, such as whether it is noise sensitive\footnote{Generalising from $10$ to $n$, the resulting family of functions is not noise sensitive. This follows from the discussion at \url{https://mathoverflow.net/questions/319439} and \url{https://mathoverflow.net/questions/324171}, together with \smartref{Theorem}{theorem_fourierclassificationofnoisesensitivity}. We do not know whether it is noise stable.}. The tools presented above do not, however, suffice, since they only apply to functions on finitely many bits. Our project will be to extend the theory to cover this case.

Most of the concepts we introduce, and many of the proofs, are directly transferred from the finite case, with various degrees of adaptation. While we do not assume the reader knows the theory in the finite case, we will often refer to it for comparison and to be clear about what is new and what is merely taken directly from that setting. Whenever we say ``the proof is as in the finite case'', or similar, the proof being referred to can be found in either \cite{garban_steif_2015} or \cite{odonnell_2014}, or both.

\section{Definition and elementary properties}

\begin{definition}
A \emph{Boolean function} is a function from $\{-1,1\}^\N$ to $\{-1,1\}$. We denote the set $\{-1,1\}$ by $\Omega$, and when convenient, we will consider $\Omega$ to be a group with multiplication, isomorphic to $\Z_2$. We will also consider functions $\Omega^\infty \to \R$.

The set $\Omega^\infty$ has an obvious uniform probability measure arising as the product measure of the uniform measure on $\Omega$ (or as the Haar measure on this compact abelian group), so we can and consistently will see $\omega\in\Omega^\infty$ as a random variable, and talk about functions on $\Omega^\infty$ lying in various $L^p$ spaces.

For a fixed function $f$ and an input $\omega = (\omega_1, \omega_2,\ldots) \in \Omega^\infty$, a set $W\subset\N$ is said to be a witness set for $\omega$ if $f(\omega)$ is almost surely determined by the bits in $W$, that is, if there exists an $A\subseteq\Omega^\infty$ such that $\Prob{\omega \in A} = 1$ and 
$$\forall \Tilde{\omega} \in A: \left(\forall i \in W: \Tilde{\omega}_i = \omega_i\right)\Rightarrow\left(f(\Tilde{\omega})=f(\omega)\right)$$

A usual Boolean function of finitely many bits thus has $[n] \equiv \{1,2,\ldots,n\}$ as witness on all inputs. Such functions are said to be finitely supported.
\end{definition}

\begin{definition}
We define the \emph{pivotal set} of a function $f$ on an input $\omega$ as being the set of all $i\in \N$ such that flipping bit $i$ changes the value of $f$ on $\omega$. That is,
$$\Pset(f)(\omega) = \left\{i\in \N: f(\omega_1,\ldots,\omega_{i-1},-\omega_i,\omega_{i+1},\ldots) \neq f(\omega)\right\}$$
A bit $i$ is said to be pivotal (on input $\omega$) if it lies in the pivotal set. The \emph{influence} of bit $k$ on a function $f$ is defined as $I_k(f) = \Prob{k \in \Pset(f)}$, and the \emph{total influence} of a function $f$ is defined as
$$I(f) = \sum_k I_k(f) = \E{\abs{\Pset(f)}}$$
\end{definition}

\begin{remark}
If $W$ is a witness set for $f$ on $\omega$, then we have $\Pset(f)(\omega) \subseteq W$, as is easily seen. In the opposite direction, we have no bounds -- it can be the case, like for the majority functions, that the pivotal set is empty but the witness set is very large. Philosophically, this shows that one should not think of the pivotal set as ``the set of bits on which $f$ depends''. 

If someone tells you what the pivotal set is, and the value of the bits in it, you know that changing the value of \emph{a single} bit outside it won't change the value of the function, but you don't necessarily know the value of the function. If you know a witness set and the value of the bits in it, however, you know that you can change \emph{any number} of bits outside the set without changing the value of the function, and you know the value of the function.
\end{remark}

\begin{definition}\label{definition_finitary_function}
A \emph{finitary Boolean function} is a function $f: \{-1,1\}^\N \to \{-1,1\}$ such that there almost surely, for a uniformly random input, exists a finite witness $W\subset \N$. For a finitary function, we make this $W$ into a well-defined random variable by choosing the least witness set according to the following ordering:
\begin{enumerate}
    \item If $\max W < \max W'$, then $W < W'$
    \item If $\max W = \max W'$ and $\abs{W} < \abs{W'}$, then $W < W'$
    \item If two witness sets have the same maximum element and the same cardinality, choose lexicographically.
\end{enumerate}

We further say that $f$ is $p$-\emph{knowable} if $\E{(\max W)^p} < \infty$. 
\end{definition}

\begin{remark}
The property of being finitary should be thought of being a function that we can, in reality, compute. It means that there exists an algorithm for the function that will halt on almost every input. It also means that we can think of the number of input bits as being only potentially infinite, not an actual infinity that we need to fit in our computer to know the value of the function.

It is thus, despite its definition seeming perhaps a bit technical, a very natural choice of class of functions to study.
\end{remark}

\begin{prop}\label{prop_pivotal_set_finite}
Any finitary Boolean function almost surely has only finitely many pivotal bits. If in addition $f$ is 1-knowable, then the total influence of $f$ is finite, that is, 
$$\sum_k I_k(f) < \infty$$
\begin{proof}
Any witness set necessarily contains all pivotal bits, and there almost surely exists a finite witness set. Thus the set of pivotal bits is almost surely finite. If additionally $\E{\max W} < \infty$, then the pivotal set is contained in a set whose expected size is finite, and thus its expected size -- which is exactly the total influence -- is also finite.
\end{proof}
\end{prop}

\subsection{Fourier analysis in this setting}
Recall that there exists a general theory of Fourier analysis on locally compact groups, of which the Fourier analysis on $\Z_2^n$ is a special case. In general, we get our Fourier analysis by noting that the characters of the group form an orthonormal basis for $L^2$ of the group. For $\Z_2$, the characters are just the functions $\chi_\emptyset: 0\mapsto 1, 1\mapsto 1$ and $\chi_{\{1\}}: 0\mapsto 1, 1 \mapsto -1$. For $\Z_2^n$, we get the characters by taking products of characters of $\Z_2$ in each coordinate.

The exact same logic prevails also for $\Z_2^\infty \cong \Omega^\infty$, noting that we take finite products of the characters. Therefore, we only state some basic results about the Fourier analysis in this setting, noting that it follows from more general theory.

\begin{prop}
The functions $\chi_S(\omega) = \prod_{i\in S} \omega_i$, indexed by finite subsets $S$ of $\N$, are the characters of the compact abelian group $\Z_2^\infty$, and thus form an orthonormal basis for $L^2(\Omega^\infty)$. In particular, for any $f \in L^2(\Omega^\infty)$, we can write
$$f = \sum_{S\subset\N} \widehat{f}(S) \chi_S$$
where $\widehat{f}(S) = \E{f\chi_S}$, and the sum runs over all \emph{finite} subsets of $\N$. We will use the convention that a sum over $S\subset\N$ implicitly runs over only finite subsets throughout. Note that this makes the sum have countably many terms, so it is an honest summation.

We also have Plancherel's and Parseval's theorems, so that for all $f, g \in L^2$
$$\norm{f}_2^2 = \sum_{S\subset\N}\widehat{f}(S)^2$$
and
$$\E{fg} = \sum_{S\subset\N} \widehat{f}(S)\widehat{g}(S)$$
\end{prop}

\subsection{Separating p-knowable, finitary, and general functions}

When defining different classes of functions, one of course wishes to know that these classes are genuinely distinct. This is of particular interest here, since the class of all functions $\{-1,1\}^\N \to \R$ is of course the entire class of random variables on $\R$, and we couldn't possibly hope to say something interesting in general about such a broad collection. So, we give results showing that general, finitary, and $p$-knowable functions are genuinely distinct classes. We also give a lemma that shows that another tentative class of functions would in fact be the entire collection of functions.

First off, not all functions are finitary. In fact, not all functions of finite total influence are finitary, so \smartref{Proposition}{prop_pivotal_set_finite} cannot be strengthened to an if and only if.

\begin{example}
Consider the function $f$ defined by that $f(\omega)=1$ if $\omega_1=1$, or if $\omega_2=\omega_3=1$, or if $\omega_4=\omega_5=\omega_6=1$, and so on, and otherwise $f(\omega) = -1$. That this function is nondegenerate is easily verified. Now, one can easily verify that it has finite total influence, and of course, for this function, there only exists a finite witness for it if $f(\omega) = 1$ -- if $f(\omega) = -1$, we need to see a $-1$ in each of infinitely many blocks. This shows that it is not a finitary Boolean function, despite having finite total influence.
\end{example}

Secondly, not all finitary functions are $p$-knowable for some $p>0$.

\begin{example}
To begin with, take $f$ to be the function that is $\omega_2$ if $\omega_1 = 1$, otherwise it is $\omega_4$ if $\omega_3=1$, and so on. That is,
$$f(\omega) = \omega_2\ind{\omega_1=1} + \ind{\omega_1=-1}\left(\omega_4\ind{\omega_3=1} + \ind{\omega_3=-1}\left(\omega_6\ind{\omega_5=1} + \ind{\omega_5=-1}\left(\cdots\right)\right)\right)$$

Since almost surely $\omega_{2n-1}=1$ for some $n$, this function is finitary, since we can then take the witness set to be all the odd-numbered bits up to $2n-1$ and bit $2n$.

Now, we wish to modify this function so that it is not $p$-knowable for any $p>0$. To do this, we replace bit $2n-1$ with a block of $a_n$ bits, and replace the condition that bit $2n-1$ equal one with the condition that the product of the bits in the block we replaced it with be one. As can be easily seen, this operation of ``splicing in extra bits'' preserves many properties of the original function.

Specifically, the new function -- $g$, say -- will still be finitary, for precisely the same reason as the original was. Now, one can easily compute that, for the original function, the probability that bit $2n$ is the final bit in the witness set is precisely $2^{-n}$. Using this, it is easy to compute that
\begin{equation*}
    \begin{split}
        \E{\left(\max W(g)\right)^p} &= \sum_{n=1}^\infty \left(n + \sum_{i=1}^n a_i\right)^p 2^{-n}
    \end{split}
\end{equation*}
and it is easy to see that we can make this sum diverge for all $p>0$ by just taking $a_n$ to grow sufficiently quickly.
\end{example}

However, it \emph{is} true that all function from $\Omega^\infty$ are approximable by functions of finite support, that is, that the set of finitely supported Boolean functions is dense in $L^1$. Thus if we want to study ``functions that are nearly functions of finitely many bits'' we cannot mean ``nearly'' in an $L^1$ sense. This follows from the following known lemma:

\begin{lemma}\label{lemma_approximation_for_nonfinitary}
Suppose $\xi_1, \xi_2, \ldots$ are random variables on some probability space $(\Omega,\F,\mu)$. Let $A\in \sigma(\xi_1, \xi_2, \ldots)$. Then, for every $\epsilon>0$, there exists an $N$ and an event $B\in\sigma(\xi_1, \ldots, \xi_N)$ so that $\mu(A\symdiff B) < \epsilon$.
\begin{proof}
Follows easily from an application of Dynkin's $\pi$-$\lambda$-lemma, or can be found in your favourite probability text.
\end{proof}
\end{lemma}

So, using that we can write a Boolean $f$ as $2\ind{f(\omega)=1}-1$ and applying our lemma to the event $f(\omega) = 1$, we get that:

\begin{cor}
Let $f:\Omega^\infty \to \Omega$ be a Boolean function. Then, for every $\epsilon>0$, there exists $N\in\N$ and a Boolean function $g: \Omega^N\to \Omega$ such that $\Prob{f \neq g} < \epsilon$.
\end{cor}

%% file: 2_noise_sensitivity_basics.tex
\section{Noise sensitivity, basics}

Having set up what our functions are, and noted that we still have a Fourier analysis, we can now move on to the main theme -- results on noise sensitivity for these functions. So, to begin with, we give a definition and some elementary results showing that this notion still behaves as one would expect from the finite case.

\begin{definition}
A sequence $f_n$ of Boolean functions is said to be noise-sensitive if $f_n(\omega)$ and $f_n(\omega^\epsilon)$ are asymptotically independent for every $\epsilon$, where $\omega$ is uniform on $\Omega^\infty$ and $\omega^\epsilon$ is $\omega$, except each bit is rerandomised with probability $\epsilon$.\footnote{Formally, this means we take a new independently sampled $\Tilde{\omega}$, and an iid sequence $X_1, X_2, \ldots$ where each $X_i$ is $1$ with probability $1-\epsilon$, and $0$ with probability $\epsilon$. We then define for each $i$
$$\omega^\epsilon_i = X_i\omega_i + (1-X_i)\Tilde{\omega_i}$$
That this matches the intuitive description of independently rerandomising bits should be clear.

We could also say that each bit is flipped independently with probability $\epsilon/2$. That is, we take an iid sequence $Y_1, Y_2,\ldots$ where $Y_i = 1$ w.p. $1-\epsilon/2$ and $Y_i = -1$ w.p. $\epsilon$, and define
$$\omega^\epsilon_i = Y_i\omega_i$$

The reader easily verifies that these two definitions are equivalent.} That is, it is noise sensitive if for every $\epsilon>0$,
$$\lim_{n\to\infty} \E{f_n(\omega)f_n(\omega^\epsilon)} - \E{f_n(\omega)}^2 = 0$$

It is said to be noise-stable if
$$\lim_{\epsilon\to 0} \sup_n \Prob{f_n(\omega) \neq f_n(\omega^\epsilon)} = 0$$
\end{definition}

\begin{prop}\label{prop_singlefunctionnoisestable}
A single function is noise stable. That is, if $f_n = f$ for all $n$, the sequence $f_n$ is noise stable.
\begin{proof}
We give a proof in the finitary case, since we can then give a more elementary proof, and defer the general case for later.

Take $f$ to be our finitary function. We wish to show that
$$\lim_{\epsilon\to 0} \Prob{f(\omega)\neq f(\omega^\epsilon)} = 0$$

So, suppose we choose $\omega$ uniformly at random. Then, by assumption, there almost surely exists a finite witness set $W_\omega$. Now, when applying an $\epsilon$-noise, if none of the bits in $W_\omega$ flip, the value of the function remains unchanged. Now, the number of bits flipped -- call it $X$ -- is conditionally binomial with parameters $\abs{W}$ and $\epsilon$. So we get that
\begin{equation*}
    \begin{split}
        \Prob{f(\omega) = f(\omega^\epsilon)} &\geq \Prob{X = 0}\\
        &= \sum_{k=0}^\infty (1-\epsilon)^k\Prob{W=k}
    \end{split}
\end{equation*}
and the limit of this expression as $\epsilon\to 0$ is $1$, by your favourite convergence theorem.
\end{proof}
\end{prop}

The above proof of course generalises to show that a sequence $f_n$ is noise-stable if the sequence $(\max W_n)_{n=1}^\infty$ is tight. This tightness is, however, not a necessary condition, as the majority functions show.

That the two concepts are called noise sensitivity and stability implies that these should be mutually exclusive concepts. This is almost true -- a constant function is trivially both noise sensitive and noise stable, and so any sequence that goes towards being constant will be both noise sensitive and noise stable. This is, fortunately, the only intersection of the two classes.

\begin{exercise}\label{exercise_stable_and_sensitive_implies_degenerate}
A sequence of functions $f_n$ is both noise sensitive and noise stable if and only if it is degenerate, in the sense that $\Var{f_n} \to 0$.
\end{exercise}

\begin{exercise}\label{exercise_neither_stable_nor_sensitive}
There are functions which are neither noise sensitive nor noise stable.
\end{exercise}

Nicely enough, since the Fourier analysis is basically the same as in the finite case, we also get the same Fourier classification of noise sensitivity and stability.

\begin{definition}
The \emph{energy at level $i$} of a function $f:\Omega^\infty \to \R$ is given by
$$E_f(i) = \sum_{S\subset\N, \abs{S} = i} \widehat{f}(S)^2$$
\end{definition}

\begin{theorem}\label{theorem_fourierclassificationofnoisesensitivity}
A sequence of Boolean functions $f_n$ is noise sensitive if and only if, for every $k>0$, the energy at non-zero levels below $k$ goes to zero with $n$.\footnote{Recall that ``level'' here refers to \emph{size} of the set of the Fourier coefficient, not to how far into the bit-sequence it lies. So a function can have the bits it is determined by run off to infinity without having the energy at finite levels going to zero.} That is, it is noise-sensitive, if and only if for every $k>0$,
$$\sum_{1\leq i\leq k} E_{f_n}(i) = \sum_{1\leq i\leq k}\sum_{S\subset\N, \abs{S}=i} \widehat{f_n}(S)^2 \to 0$$
as $n\to\infty$.
\begin{proof}
To begin with, note that
\begin{equation*}
    \begin{split}
        \E{f_n(\omega)f_n(\omega^\epsilon)} &= \E{\left(\sum_{S_1 \subset \N} \widehat{f_n}(S_1)\chi_{S_1}(\omega)\right)\left(\sum_{S_2 \subset \N} \widehat{f_n}(S_2)\chi_{S_2}(\omega^\epsilon)\right)}\\
        &= \sum_{S\subset\N} \widehat{f_n}(S)^2\E{\chi_S(\omega)\chi_S(\omega^\epsilon)}\\
        &= \sum_{S\subset \N} \widehat{f_n}(S)^2(1-\epsilon)^{\abs{S}}
    \end{split}
\end{equation*}
and of course $\E{f_n(\omega)} = \widehat{f}_n(\emptyset)$, so that
\begin{equation*}
    \begin{split}
        \lim_{n\to\infty} \E{f_n(\omega)f_n(\omega^\epsilon)} - \E{f_n(\omega)}^2 &= \lim_{n\to\infty}\left(\sum_{S\subset \N}\widehat{f_n}(S)^2(1-\epsilon)^{\abs{S}}\right) - \widehat{f_n}(\emptyset)^2\\
        &= \lim_{n\to\infty}\sum_{S\subset \N, S \neq \emptyset}\widehat{f_n}(S)^2(1-\epsilon)^{\abs{S}}\\
        &= \lim_{n\to\infty}\sum_{k=1}^\infty \sum_{S\subset \N, \abs{S}=k} \widehat{f_n}(S)^2(1-\epsilon)^k
    \end{split}
\end{equation*}
and that this expression is zero iff our hypothesis holds is easily seen. 
\end{proof}
\end{theorem}

\begin{theorem}\label{theorem_fourierclassificationofnoisestability}
A sequence of Boolean functions $f_n$ is noise stable if and only if, for every $\epsilon>0$, there exists a $k$ such that for all $n$,
$$\sum_{i=k}^\infty \sum_{\abs{S}=i} \widehat{f_n}(S)^2 < \epsilon$$
\begin{proof}
The proof of this is exactly the same as in the finite case, so since it gets very messy and we will never use this result, the proof is omitted.
\end{proof}
\end{theorem}

\begin{cor}
\smartref{Proposition}{prop_singlefunctionnoisestable} holds for all functions, not just finitary ones.
\end{cor}

\begin{cor}
This gives solutions to \smartref{Exercise}{exercise_stable_and_sensitive_implies_degenerate} and \smartref{Exercise}{exercise_neither_stable_nor_sensitive}. They can also be solved by direct calculation, without resorting to Fourier analysis.
\end{cor}

%% file: 3_0_finite_approximations_and_hypercontractivity.tex
\section{Approximations of finitary and non-finitary functions, and hypercontractivity}

If a sequence $f_n$ is noise sensitive, this can be phrased as that $f_n$ and $f_n$ with noise applied are asymptotically independent. It thus makes sense to define more precisely what we mean by a noised version of a function, which can be done more generally than we have defined noise sensitivity.

\begin{definition}
We define, for each $\rho\in [0,1]$ and for all $p\in[1,\infty)$, the \emph{noise operator} $T_\rho: L^p(\Omega^\infty)\to L^p(\Omega^\infty)$ by
$$(T_\rho f)(\omega) = \E{f(\omega^{1-\rho}) \given \omega}$$
Whenever we apply $T_\rho$ to a function from $\Omega^n$, we mean the obvious restriction of $T_\rho$ induced by seeing $L^p(\Omega^n)$ as a subspace of $L^p(\Omega^\infty)$.

It is easily seen that these operators are linear. That they indeed map into $L^p$ will be seen later.
\end{definition}

As can be easily seen, this gives $\E{f(\omega)f(\omega^\epsilon)} = \inpr{f}{T_{1-\epsilon} f}$. Now, the question is, how should we expect the noised version of $f$ to behave? Conditional on $\omega$, the distribution of $\omega^\epsilon$ will be centered around $\omega$, mostly living on points close to $\omega$ in the Hamming distance. Noising can thus be seen as a local averaging operation, and should thus act in a regularising and smoothing way on a function -- its extrema are replaced by local averages, which of course by definition are no larger than the extremum itself.

One can also think of the noise as being the product of a random walk on the Hamming cube -- we start at $\omega$, and randomly take steps along the edges for some time corresponding to the amount of noise. So the noising can also be thought of as a type of diffusion, which should again act in a regularising way.

It turns out that this is indeed the case, and we get a very useful bound of the $L^2$-norm of the noised function by the $L^{1+\rho^2}$-norm of the original function. In the finite setting, the statement is as follows:\footnote{For a proof, and discussion of whom exactly the result should be attributed to, see \cite{garban_steif_2015}.}

\begin{theorem}[Hypercontractivity, finite case]
For any $n\in\N$ and $\rho\in[0,1]$, it holds for every $f: \Omega^n \to \R$ that
$$\norm{T_\rho f}_2 \leq \norm{f}_{1+\rho^2}$$
\end{theorem}

Of course, in the finite setting, a function from $\Omega^n$ can only take at most $2^n$ different values, and so all functions are in $L^\infty$. In the case of functions from $\Omega^\infty$, this is no longer the case, so a regularising bound like this becomes potentially even more interesting, in that it can ``promote'' a function from $L^q$ to $L^2$ for $q<2$.

In order to extend this result from the finite case, we will establish various results on how a function from $\Omega^\infty$ can be approximated by one from $\Omega^n$. Once we have good enough convergence results for these approximations, and have seen that they interact very well with the noise operator, extending hypercontractivity will be straightforward. These approximation results will keep reappearing in later sections as well, since they allow us to reduce questions about functions on infinitely many bits to questions about functions on finitely many bits.

\input{3_1_approximations.tex}
\input{3_2_hypercontractivity.tex}
\input{3_3_apps_and_noise_sensitivity.tex}

%% file: 3_1_approximations.tex
\subsection{Approximations of Boolean functions}
\begin{definition}\label{definition_fourier_projection_approximation}
We define the \emph{approximation operator} $A_n: L^1(\Omega^\infty) \to L^1(\Omega^n)$ by
$$(A_n f)(\omega) = \E{f \given \omega_1,\omega_2,\ldots,\omega_n}$$

So $A_n$ is the projection onto the space of functions depending only on the first $n$ bits. If $f\in L^2$, this gives a simple formula for $A_n$: Since the $\chi_S$ for $S\subseteq[n]$ are an orthonormal basis for the subspace we're projecting onto, we get that
$$A_n f = \sum_{S\subseteq[n]}\widehat{f}(S)\chi_S$$
\end{definition}

\begin{definition}\label{definition_boolean_approximation}
We also define the \emph{Boolean approximation operator} $B_n$ taking finitary Boolean functions to Boolean functions on $n$ bits by
$$B_n f = \ind{W\subseteq[n]}f - \ind{W\not\subseteq[n]}$$
That is, $B_n f(\omega) = f(\omega)$ if we can know the value of $f$ on $\omega$ using only the first $n$ bits, otherwise it is $-1$.

Note that this operator is \emph{not} linear, unlike $A_n$.
\end{definition}

\begin{remark}
In general, $A_n f$ will not itself be a Boolean function! So it is very nice to work with on the Fourier side, but we lose all the combinatorial structure of the function, so we cannot easily reason about it in any non-analytical way. Similarly, $B_n f$ will be a Boolean function that inherits most of its combinatorial structure from $f$, but there is no good connection between its Fourier coefficients and those of $f$. So one has to choose on which side one wants to be approximating, since we can't get good properties on both sides.
\end{remark}

Getting convergence for the Boolean approximations is straightforward, since we can use the strength of our assumption that $f$ is finitary.

\begin{prop}\label{prop_boolean_approximations_converge}
Suppose $f$ is a finitary Boolean function. Then $B_n f \to f$ in $L^p$ for every $p\in[1,\infty)$.
\begin{proof}
Since $W$ is almost surely finite, the event $W\subseteq[n]$ increases towards being the entire probability space, and so we have $B_n f \nearrow f$ pointwise. Since $B_nf, f\in L^p$ the monotone convergence theorem now gives the result.
\end{proof}
\end{prop}

In order to get that $A_n f$ converges to $f$ even when $f$ is not finitary, we will need some heavier analytical artillery.

\begin{theorem}[Doob's Martingale Convergence Theorem in $L^p$]\label{theorem_doob_martingale_convergence}
Suppose $(\Omega,\F,\mathbb{P})$ is some probability space, and $X_k$ is a martingale. For all $p\in[1,\infty)$, if there exists a constant $K$ such that $\norm{X_k}_p \leq K$ for all $k$, then there exists a random variable $X \in L^p$ such that $X_k \to X$ almost surely and in $L^p$. 
\end{theorem}

\begin{theorem}[Levy's Upward Theorem]\label{theorem_levy_upward}
Suppose $(\Omega,\F,\mathbb{P})$ is some probability space, and $\F_n \nearrow \F$ is some filtration. Then for any $X\in L^1$, it holds that $\E{X \given \F_n} \to X$ almost surely.
\end{theorem}

\begin{prop}\label{prop_approximations_converge}
For every $p\in[1,\infty)$, it holds for every $f\in L^p(\Omega^\infty)$ that $A_nf \to f$ almost surely and in $L^p$.
\begin{proof}
Take $\F_n$ to be the sigma algebra generated by the first $n$ coordinates. We then get that $A_n f = \E{X \given \F_n}$, and so \smartref{Theorem}{theorem_levy_upward} gives us the almost sure convergence.

For the convergence in $L^p$, note that Jensen's inequality gives us that
$$\norm{A_n f}_p^p = \E{\abs{\E{f \given \F_n}}^p} \leq \E{\E{\abs{f}^p \given \F_n}} = \norm{f}_p^p$$
so we can take $K=\norm{f}_p$ in \smartref{Theorem}{theorem_doob_martingale_convergence}, and so we get convergence in $L^p$ as well.
\end{proof}
\end{prop}

\begin{exercise}
For a finitary Boolean function $f: \Omega^\infty \to \Omega$, we do not need any heavy machinery, and can prove in a few lines that $A_n f \to f$ in $L^p$ for all $p$ by using that $A_n f = f$ conditional on $\max W \leq n$.
\end{exercise}

%% file: 3_2_hypercontractivity.tex
\subsection{Hypercontractivity}

Having set up our approximations, we can now continue to describe how the noise operator works, in order to get that it commutes with our projection approximation. This in turn will make extending hypercontractivity from the finite to the finitary case straightforward.


\begin{lemma}\label{lemma_noise_short_map}
For all $p\in[1,\infty)$ and all $\rho$, we have for all $f\in L^p(\Omega^\infty)$ that
$$\norm{T_\rho f}_p \leq \norm{f}_p$$
and in particular, $T_\rho$ is a bounded linear operator from $L^p$ to $L^p$.
\begin{proof}
That it is linear is immediate by linearity of expectation. For the boundedness, we have by definition that
\begin{equation*}
    \begin{split}
        \norm{T_\rho f}_p^p &= \E{\abs{T_\rho f}^p}\\
        &= \E{\abs{\E{f(\omega^{1-\rho}) \given \omega}}^p}
    \end{split}
\end{equation*}
and Jensen's inequality then gives
$$\E{\abs{\E{f(\omega^{1-\rho}) \given \omega}}^p} \leq \E{\E{\abs{f(\omega^{1-\rho})}^p \given \omega}} = \norm{f}_p^p$$
where we used that $\omega$ and $\omega^{1-\rho}$ have the same distribution in the last equality.
\end{proof}
\end{lemma}

The following is not immediately necessary to prove hypercontractivity, but will be very useful when applying it, since it will be applied by massaging a Fourier-style expression into being a noising of something, and then applying hypercontractivity to get a better norm bound on said expression. It also even further justifies seeing noising as a type of diffusion -- it acts exactly like heat diffusion in damping the energy at each level proportional to its ``frequency''.

\begin{lemma}\label{lemma_spectral_formula_for_noise_operator}
Assume $f\in L^2(\Omega^\infty)$. $T_\rho f$ is also given by
$$T_\rho f = \sum_{S\subset \N} \rho^{\abs{S}}\widehat{f}(S)\chi_S$$
\begin{proof}
Since $T_\rho$ is, by \smartref{Lemma}{lemma_noise_short_map}, a continuous linear operator, we can write
$$T_\rho f = T_\rho\left(\sum_{S\subset\N} \widehat{f}(S)\chi_S\right) = \sum_{S\subset\N} \widehat{f}(S)T_\rho\left(\chi_S\right)$$
so it suffices to check how $T_\rho$ acts on the characters.

So, fix some $S\subset\N$. We then have
$$T_\rho\chi_S(\omega) = \E{\chi_S(\omega^{1-\rho}) \given \omega}$$
and, letting $F$ be the set of bits rerandomised in $S$, we see
\begin{equation*}
    \begin{split}
        \E{\chi_S(\omega^{1-\rho}) \given \omega} &= \E{\chi_{S\setminus F}(\omega)\chi_F(\omega^{1-\rho}) \given \omega}\\
        &= \E{\chi_{S\setminus F}(\omega)\chi_F(\omega^{1-\rho}) \given \omega, F\neq\emptyset}\Prob{F\neq \emptyset \given \omega}\\
        &\qquad + \E{\chi_{S\setminus F}(\omega)\chi_F(\omega^{1-\rho}) \given \omega, F=\emptyset}\Prob{F = \emptyset \given \omega}
    \end{split}
\end{equation*}
where the second term in this sum is precisely $\rho^{\abs{S}}\chi_S(\omega)$, and the first term is easily seen to be zero by independence.
\end{proof}
\end{lemma}

\begin{lemma}\label{lemma_approximation_commutes_with_noise}
For all $\rho\in[0,1]$ and all $p\in[1,\infty)$, $A_n$ and $T_\rho$ commute as operators on $L^p$. That is, for all $f\in L^p(\Omega^\infty)$,
$$A_nT_\rho f = T_\rho A_n f$$
\begin{proof}
For $p\geq 2$, we can immediately read this off from the Fourier representations of $A_n$ and $T_\rho$. For $p<2$, we no longer have such a representation, and instead need to look at their definitions in terms of conditional expectation.

So, to begin with, let $\epsilon = 1-\rho$, and let $\xi_i$ be iid, with $\xi_i = -1$ w.p. $\frac{\epsilon}{2}$ and $1$ otherwise, so that we can let $\omega^{1-\rho}=\omega^\epsilon=\omega\cdot\xi$. We also denote $\omega_1, \omega_2, \ldots, \omega_n$ by $\omega_{[n]}$. With these notational conveniences established, we can now calculate
$$A_nT_\rho f = \E{T_\rho f \given \omega_{[n]}} = \E{\E{f(\omega\cdot\xi) \given \omega} \given \omega_{[n]}} = \E{f(\omega\cdot\xi) \given \omega_{[n]}}$$
and
\begin{equation*}
    \begin{split}
        T_\rho A_nf &= \E{A_nf((\omega\cdot\xi)_{[n]}) \given \omega}\\
        &= \E{A_nf((\omega\cdot\xi)_{[n]}) \given \omega_{[n]}}\\
        &= \E{\E{f(\omega\cdot\xi) \given (\omega\cdot\xi)_{[n]}}\given\omega_{[n]}}\\
        &= \E{\E{f(\omega\cdot\xi) \given \omega_{[n]}, \xi_{[n]}}\given\omega_{[n]}}\\
        &= \E{f(\omega\cdot\xi) \given \omega_{[n]}}
    \end{split}
\end{equation*}
where the second equality follows from that $A_n f$ is independent of bits after bit $n$, and the fourth equality follows from that learning which bits were flipped and what the original states were gives the exact information as learning the final states of the bits for computing $f$ on the rerandomised sequence. So they are equal, as desired.
\end{proof}
\end{lemma}

\begin{theorem}[Hypercontractivity]\label{theorem_hypercontractivity}
For every $\rho\in[0,1]$, we have for any function $f\in L^{1+\rho^2}(\Omega^\infty)$ that
$$\norm{T_\rho f}_2 \leq \norm{f}_{1+\rho^2}$$
and so in particular $T_\rho$ is a weak contraction from $L^{1+\rho^2}$ into $L^2$.
\begin{proof}
The finite case of hypercontractivity gives us for every $n\in\N$ that
$$\norm{T_\rho A_n f}_2 \leq \norm{A_n f}_{1+\rho^2}$$
\smartref{Lemma}{lemma_approximation_commutes_with_noise} lets us swap the order of noising and approximation, giving
\begin{equation}\label{eq_finite_case_hypercont_proof}
    \norm{A_n\left( T_\rho f\right)}_2 \leq \norm{A_n f}_{1+\rho^2}
\end{equation}

Now, \smartref{Proposition}{prop_approximations_converge} tells us that the right hand side of this converges, but we cannot immediately apply it to the left hand side, since we do not know that $T_\rho f \in L^2$. However, that $A_n f \to f$ in $L^{1+\rho^2}$ implies in particular that $\sup_n \norm{A_n f}_{1+\rho^2} < \infty$, and so we can take this supremum as the $K$ in \smartref{Theorem}{theorem_doob_martingale_convergence}. Thus we can take the limit as $n\to\infty$ on both sides of this inequality, getting the desired inequality in the infinite case as well.
\end{proof}
\end{theorem}

%% file: 3_3_apps_and_noise_sensitivity.tex
\subsection{Approximations and noise sensitivity}

We now give a result that essentially states that noise sensitivity is, for finitary functions, something that happens in the first finitely many bits. This will not be immediately useful, but will appear as a useful tool again later. First, however, we need a small lemma:

\begin{lemma}\label{lemma_noise_autocorr_boolean_converge}
Suppose $f$ is a finitary function. Then
$$\E{B_m f(\omega) B_m f(\omega^\epsilon)} \to \E{f(\omega)f(\omega^\epsilon)}\qquad\text{as }m\to\infty$$
\begin{proof}
Let
$$K_m = \{W(f)(\omega) \subseteq[m], W(f)(\omega^\epsilon) \subseteq[m]\}$$
It is easily seen that $\Prob{K_m} \to 1$. So, we can calculate
\begin{equation*}
    \begin{split}
        \E{B_mf(\omega)B_mf(\omega^\epsilon)} &= \E{B_mf(\omega)B_mf(\omega^\epsilon)\ind{K_m}} + \E{B_mf(\omega)B_mf(\omega^\epsilon)\ind{K_m^c}}\\
        &= \E{f(\omega)f(\omega^\epsilon)\ind{K_m}} + \E{B_mf(\omega)B_mf(\omega^\epsilon)\ind{K_m^c}}\\
        &\xrightarrow{m\to\infty} \E{f(\omega)f(\omega^\epsilon)}
    \end{split}
\end{equation*}
which is what was to be shown.
\end{proof}
\end{lemma}

\begin{prop}\label{prop_approximations_noise_sensitive_implies_noise_sensitive}
Suppose $f_n$ is some sequence of finitary Boolean functions. Then the following are equivalent:
\begin{enumerate}
    \item For every sequence of integers $m_n$ going to infinity sufficiently quickly, the sequence $B_{m_n}f_n$ is noise sensitive. 
    
    That is, there exists some sequence $r_n\to\infty$ such that for any sequence $m_n$ that goes to infinity at least as fast as $r_n$ (that is, for all $n$, $m_n \geq r_n$), the sequence $B_{m_n} f_n$ is noise sensitive.
    \item $f_n$ is noise sensitive.
\end{enumerate}
\begin{proof}
To cut down on the volume of formulas, for each function $f$, let
$$\Xi(f) = \E{f(\omega)f(\omega^\epsilon)} - \E{f(\omega)}^2$$
so that a sequence $f_n$ is noise sensitive iff $\Xi(f_n) \to 0$.

We have from \smartref{Proposition}{prop_boolean_approximations_converge} and \smartref{Lemma}{lemma_noise_autocorr_boolean_converge} that for each $n$
\begin{equation*}
    \Xi(B_m f_n) \xrightarrow{m\to\infty} \Xi(f_n)
\end{equation*}

From this we see that we can for each $n$ choose an $a_n$ such that for all $m\geq a_n$
$$\abs{\Xi(B_m f_n) - \Xi(f_n)} < \frac{1}{n}$$

So, suppose that there exists some sequence $r_n$ so that for any sequence $m_n \geq r_n$ the sequence $B_{m_n} f_n$ is noise sensitive. Letting $b_n = \max(r_n, a_n)$, we can compute
$$\abs{\Xi(f_n)} \leq \abs{\Xi(f_n) - \Xi(B_{b_n} f_n} + \abs{\Xi(B_{b_n} f_n)} < \frac{1}{n} + \abs{\Xi(B_{b_n} f_n)} \to 0$$
where we used that $B_{b_n} f$ is noise sensitive since $b_n \geq r_n$ to get that $\Xi(B_{b_n} f) \to 0$. Thus the sequence $f_n$ is noise sensitive.

Now suppose $f_n$ is noise sensitive, and let $b_n$ be some sequence such that $b_n \geq a_n$ for all $n$. We can then compute
$$\abs{\Xi(B_{b_n} f_n)} \leq \abs{\Xi(B_{b_n} f_n) - \Xi(f_n)} + \abs{\Xi(f_n)} < \frac{1}{n} + \abs{\Xi(f_n)} \to 0$$
where we used the assumption that $f_n$ is noise sensitive to get that $\Xi(f_n) \to 0$. Thus the sequence $B_{b_n} f_n$ is noise sensitive.
\end{proof}
\end{prop}

%% file: 4_0_influences_and_BKS.tex
\section{Influences, the BKS noise sensitivity theorem, and the KKL maximal influence theorems}

Having proven hypercontractivity, we can now move on to proving our versions of the BKS noise sensitivity theorem for finitary functions. To do this, we will first need to say more about how influences work for functions on infinitely many bits, since we will need approximation results for these too in order to transfer noise sensitivity results from the finite case.

\input{4_1_influences_and_approximations.tex}

\input{4_2_BKS.tex}
\input{4_3_KKL.tex}

%% file: 4_1_influences_and_approximations.tex
\subsection{Influences and how they interact with approximations}

To begin with, we need to establish that the Fourier formula for the influences of a function from the finite case also holds in our setting.

As a notational convenience, we use $\omega^k$ to mean $\omega$ flipped at bit $k$, that is,
$$\omega^k = (\omega_1, \omega_2, \ldots, \omega_{k-1}, -\omega_k, \omega_{k+1}, \ldots)$$
No confusion with the notation $\omega^\epsilon$ for a rerandomised $\omega$ should come from this, since they don't occur together, and it should be clear when we have an integer and when we have a small $\epsilon$.

\begin{lemma}\label{lemma_fourier_formula_for_influences}
If $f$ is any Boolean function, not necessarily finitary, it holds that
$$I_k(f) = \sum_{S \ni k} \widehat{f}(S)^2$$
and
$$I(f) = \sum_S \abs{S}\widehat{f}(S)^2$$
\begin{proof}
We can calculate
\begin{equation*}
    \begin{split}
        I_k(f) &= \Prob{f(\omega) \neq f(\omega^k)}\\
        &= \Prob{f(\omega)f(\omega^k) = -1} = \frac{1 - \E{f(\omega)f(\omega^k)}}{2}
    \end{split}
\end{equation*}
and for the expected value we can Fourier expand $f$, and using the usual orthogonality of the characters get that
\begin{equation*}
    \begin{split}
        \E{f(\omega)f(\omega^k)} &= \E{\left(\sum_S\widehat{f}(S)\chi_S(\omega)\right)\left(\sum_S\widehat{f}(S)\chi_S(\omega^k)\right)}\\
        &= \sum_S \widehat{f}(S)^2\E{\chi_S(\omega)\chi_S(\omega^k)}\\
        &= \sum_{S\not\ni k} \widehat{f}(S)^2 - \sum_{S\ni k} \widehat{f}(S)^2\\
        &= \left(1 - \sum_{S\ni k}\widehat{f}(S)^2\right) - \sum_{S\ni k} \widehat{f}(S)^2
    \end{split}
\end{equation*}
where we in the last line used that $\sum_S \widehat{f}(S)^2 = \E{f^2} = 1$. Inserting this into our expression for $I_k(f)$ in terms of this expected value gets exactly the desired formula, and summing this over all $k$ gives the formula for $I(f)$.
\end{proof}
\end{lemma}

We also note that influences can be written in terms of a discrete derivative.
\begin{definition}
For a Boolean function $f: \Omega^\infty \to \Omega$, we define
$$\nabla_k f: \omega \mapsto \frac{f(\omega) - f(\omega^k)}{2}$$

A simple computation gives that
$$\nabla_k f(\omega) = \frac{1}{2}\sum_S \widehat{f}(S)\left(\chi_S(\omega) - \chi_S(\omega^k)\right) = \sum_{S\ni k} \widehat{f}(S)\chi_S(\omega)$$
so that
$$\widehat{\nabla_k f}(S) = \begin{cases}\widehat{f}(S) &S \ni k\\
0&\text{otherwise}\end{cases}$$

Observe that $I_k(f) = \norm{\nabla_k f}_1$, and since $\nabla_k f \in \{-1,0,1\}$, we in fact get that $I_k(f) = \norm{\nabla_k f}_p^p$ for all $p\in [1,\infty)$. Taking $p=2$ and applying Parseval's formula gives an alternative proof of \smartref{Lemma}{lemma_fourier_formula_for_influences}.
\end{definition}

\begin{definition}
For any Boolean function $f$, we define
$$H(f) = \sum_k I_k(f)^2$$

Note that this can equivalently, but more combinatorially, be defined as
$$H(f) = \E{\abs{\Pset(f)(\omega)\cap\Pset(f)(\Tilde{\omega})}}$$
where $\omega$ and $\Tilde{\omega}$ are independent.
\end{definition}

\begin{lemma}\label{lemma_knowability_bounds_influences}
Suppose $f$ is a $p$-knowable Boolean function. We then have that
$$I_k(f) \leq \frac{\E{\left(\max W\right)^p}}{k^p}$$
for every $k$.
\begin{proof}
Fix some $k$. An easy calculation using the fact that the pivotal set is always a subset of any witness set gives that
\begin{equation*}
    \begin{split}
        I_k(f) &= \Prob{k \in \Pset(f)}\\
        &\leq \Prob{k \in W}\\
        &\leq \Prob{k \leq \max W}\\
        &\leq \frac{\E{\left(\max W\right)^p}}{k^p}
    \end{split}
\end{equation*}
with the last inequality being Markov's inequality.
\end{proof}
\end{lemma}

\begin{cor}\label{cor_influence_tail_bound}
If $f$ is $p$-knowable and $q>0$ is such that $pq>1$, then it holds that
$$\sum_{k>n} I_k(f)^q \leq \E{\left(\max W\right)^p}^q\frac{n^{1-pq}}{pq-1}$$
\end{cor}

\begin{cor}
If $f$ is $p$-knowable for some $p>\frac{1}{2}$, then $H(f) < \infty$.
\end{cor}

\begin{prop}\label{prop_approximations_and_influences}
For any finitary Boolean function $f$, $I_k(B_n f) \to I_k(f)$. If additionally $f$ is $p$-knowable for some $p > \frac{1}{2}$, then $H(B_n f) \to H(f)$, and if it is $p$-knowable for some $p>1$, then $I(B_n f)\to I(f)$.
\begin{proof}
For every fixed $k$ we have
\begin{equation*}
    \begin{split}
        \left\{B_n f(\omega) \neq B_n f(\omega^k)\right\} = \big\{f(\omega) &\neq f(\omega^k), W(\omega) \subseteq[n], W(\omega^k) \subseteq[n]\big\}\\
        &\cup \left\{B_n f(\omega) \neq B_n f(\omega^k), W(\omega)\cup W(\omega^k) \not\subseteq[n]\right\}
    \end{split}
\end{equation*}
and since the witness sets are almost surely finite, the first of the two events in the union increases towards being the event $\{f(\omega) \neq f(\omega^k)\}$, and the second goes towards being a null event. So, for each fixed $k$, $I_k(B_n f) \to I_k(f)$, as desired.

Now, if we additionally assume that $f$ is $p$-knowable, it looks a whole lot like we could just use \smartref{Lemma}{lemma_knowability_bounds_influences} to get the dominated convergence theorem to give us the result. However, what we in fact need is to get that bound \emph{uniformly across all $B_n f$} as well. Of course, that $B_n f$ is $p$-knowable is trivial, so what remains to establish is that
$$\sup_n \E{\left(\max W(B_n f)\right)^p} < \infty$$

Some simple calculations give us that
\begin{equation*}
    \begin{split}
        \sup_n \E{\left(\max W(B_n f)\right)^p} &\leq \E{\sup_n \left(\max W(B_n f)\right)^p}\\
        &= \E{\left(\sup_n \max W(B_n f)\right)^p}\\
        &= \E{\left(\max \bigcup_n W(B_n f)\right)^p}
    \end{split}
\end{equation*}

Now, in fact, it holds almost surely that $\max \bigcup_n W(B_n f) = \max W(f)$. To see this, consider some fixed $\omega$. Clearly, for $n < \max W(f)(\omega)$, we must have $\max W(B_n f)(\omega) < \max W(f)(\omega)$, since the witness set for $B_n f$ is trivially contained in $[n]$. But for $n\geq \max W(f)(\omega)$, we have $W(B_n f)(\omega) = W(f)(\omega)$ by construction.  So this big union is in fact a union of a finite number of sets which lie before $\max W(f)(\omega)$, and then infinitely many copies of $W(f)(\omega)$.

So, what we have established is that
$$\sup_n \E{\left(\max W(B_n f)\right)^p} \leq \E{\left(\max W(f)\right)^p}$$
which using \smartref{Lemma}{lemma_knowability_bounds_influences} gives us that
$$I_k(B_n f) \leq \frac{\E{\left(\max W(B_n f)\right)^p}}{k^p} \leq \frac{\E{\left(\max W(f)\right)^p}}{k^p}$$
which is a bound that is independent of $n$, as desired.

If we have $p>\frac{1}{2}$, then this bound is square-summable, and we get $H(B_n f)\to H(f)$ by the dominated convergence theorem. Likewise, if $p>1$, the bound is summable and we get $I(B_n f)\to I(f)$ also by dominated convergence.
\end{proof}
\end{prop}

%% file: 4_2_BKS.tex
\subsection{An inequality relating influences and noise sensitivity, and the Benjamini-Kalai-Schramm theorem}

In the finite setting, there is a clear relation between noise sensitivity and influences, that is quantified by the following result:
\begin{theorem}[\cite{KK_13}]\label{theorem_KK_inequality}
There exist constants $a, b > 0$ such that we, for all Boolean functions $f: \Omega^n \to \Omega$ and for every $\epsilon>0$, have
$$\E{f(\omega)f(\omega^\epsilon)} - \E{f(\omega)}^2 \leq a(H(f))^{b\epsilon}$$
\end{theorem}

This immediately gives that if $H(f_n)\to 0$ for a sequence of Boolean functions $f_n: \Omega^{m_n} \to \Omega$, then $f_n$ is noise sensitive, which is the original noise sensitivity result of Benjamini, Kalai and Schramm\cite{BKS_theorem}. 

This result transfers quite directly to the finitary setting:

\begin{theorem}\label{theorem_finitary_KK_inequality}
Suppose $f$ is some Boolean function that is $p$-knowable for some $p>\frac{1}{2}$. There exist universal constants $a$ and $b$, not depending on $p$ or $f$, such that for all $\epsilon>0$
$$\E{f(\omega)f(\omega^\epsilon)} - \E{f(\omega)}^2 \leq a(H(f))^{b\epsilon}$$
\begin{proof}
We have, using the finite version of this result (\smartref{Theorem}{theorem_KK_inequality}), that for each $m\in\N$
$$\E{B_mf(\omega)B_mf(\omega^\epsilon)} - \E{B_mf(\omega)}^2 \leq a(H(B_mf))^{b\epsilon}$$

Since $p>\frac{1}{2}$, \smartref{Proposition}{prop_approximations_and_influences} gives that $H(B_m f) \to H(f)$. We also have from \smartref{Lemma}{lemma_noise_autocorr_boolean_converge} and \smartref{Proposition}{prop_boolean_approximations_converge} that the left hand side of the inequality converges. So, taking the limit as $m\to\infty$ gives the desired result.
\end{proof}
\end{theorem}

A version of the Benjamini-Kalai-Schramm noise sensitivity theorem is now an immediate corollary:

\begin{theorem}[BKS for finitary functions, qualitative version]\label{theorem_finitary_BKS_qualitative}
Let $f_n$ be some sequence of finitary functions, such that each $f_n$ is $p_n$-knowable for some $p_n > \frac{1}{2}$. Suppose
$$H(f_n) \to 0$$
Then $f_n$ is noise sensitive.
\end{theorem}

\begin{exercise}
\smartref{Theorem}{theorem_finitary_BKS_qualitative} could also be proved directly from the finite case of the same theorem, using \smartref{Proposition}{prop_approximations_and_influences} and  \smartref{Proposition}{prop_approximations_noise_sensitive_implies_noise_sensitive}.
\end{exercise}

It is natural to ask oneself whether we actually need $f_n$ to be finitary for this theorem to hold, since all the concept involved are well-defined also for non-finitary functions. We do not know of any example that shows that this hypothesis is needed, but there does not appear to be a way to extend the proof we have to the non-finitary case.

The main issue seems to be that $H$ is not a continuous functional on $L^2$ if we extend it using the Fourier formula for $I_k(f)$. In particular, one can compute that
$$H(f) = \sum_{S, S'} \abs{S\cap S'} \left(\widehat{f}(S)\widehat{f}(S')\right)^2$$
so, if we let $f_n = \frac{1}{n^{1/4}}\chi_{[n]}$, we see that $\norm{f_n}_2 = \frac{1}{n^{1/4}}$, while $H(f_n) = 1$, so $f_n\to 0$ while $H(f_n) \to 1 \neq 0 = H(0)$.

So we cannot extend our results from finitary to non-finitary by any continuity arguments, which seems to leave any direct path from the finite result to the general result closed. 

There is, in the finite setting, also a version of the BKS theorem that states that if $H(f)$ not only goes to zero but does so at an inversely polynomial rate, then the energy at levels below $\log(n)$ goes to zero with $n$, which implies noise sensitivity. 

The proof of this version of the theorem does transfer almost directly from the finite case, since the core ingredient in the proof is hypercontractivity, which we have extended from the finite setting. Unfortunately, we need to at one stage get a bound on the tail of the influences, which requires us to again assume $p$-knowability for $p>\frac{1}{2}$. Since both the hypotheses and conclusion are stronger than in \smartref{Theorem}{theorem_finitary_BKS_qualitative}, neither of the theorems implies the other.

\begin{theorem}[BKS for finitary functions, logarithmic version]\label{theorem_finitary_BKS_logarithmic}
Let $f_n$ be some sequence of $p$-knowable Boolean functions for some $p>\frac{1}{2}$, and let $\mu_n = \E{\left(\max W(f_n)\right)^p}$. Assume $\mu_n \to \infty$.

Suppose there exists some $\gamma > 0$ such that $H(f_n) \leq \mu_n^{-\gamma}$. Then we have
$$\sum_{1 \leq \abs{S} \leq  \log\left(\mu_n\right)} \widehat{f_n}(S)^2 \to 0$$
and in particular $f_n$ is noise sensitive, by \smartref{Theorem}{theorem_fourierclassificationofnoisesensitivity}.
\begin{proof}
Let $m_n$ be some sequence of integers going to infinity and $\rho \in (0,1)$ a constant, to be specified later. Then
\begin{equation*}
    \begin{split}
        \sum_{1 \leq \abs{S} \leq m_n} \hat{f_n}(S)^2 &\leq \sum_{1 \leq \abs{S} \leq m_n} \abs{S}\hat{f_n}(S)^2\\
        &= \sum_k \sum_{1 \leq \abs{S} \leq m_n} \widehat{\nabla_k f_n}(S)^2
        = \sum_k \sum_{1 \leq \abs{S} \leq m_n} \left(\frac{\rho^{\abs{S}}}{\rho^{\abs{S}}}\widehat{\nabla_k f_n}(S)\right)^2\\
        &\leq \sum_k \rho^{-2m_n} \sum_{1 \leq \abs{S} \leq m_n} \left(\rho^{\abs{S}}\widehat{\nabla_k f_n}(S)\right)^2\\
        &\leq \sum_k \rho^{-2m_n}\norm{T^\infty_\rho \nabla_k f_n}_2^2\\
        &\leq \sum_k \rho^{-2m_n}\norm{\nabla_k f_n}_{1+\rho^2}^2
    \end{split}
\end{equation*}
where we used hypercontractivity in the last inequality.

Now, since $f_n$ is Boolean, we have $\norm{\nabla_k f_n}_{1+\rho^2} = \norm{\nabla_k f_n}_2^{2/(1+\rho^2)}$. We would like to use Hölder's inequality to get a sum of squared influences free, but unfortunately, we cannot do so directly, since we are summing infinitely many terms. So, we cut the sum off at some $N_n$ -- to be chosen later -- and handle the two parts separately. Below $N_n$, we can apply Hölder's inequality, and it will go to zero because the sum of the squared influences goes to zero, while above $N_n$, all influences are necessarily small, so that term vanishes as well.
\begin{equation*}
    \begin{split}
         \sum_{1 \leq \abs{S} \leq m_n} \hat{f_n}(S)^2 &\leq \sum_k \rho^{-2m_n}\norm{\nabla_k f_n}_2^{4/(1+\rho^2)}\\
         &= \rho^{-2m_n}\sum_k \left(I_k(f_n)\right)^{2/(1+\rho^2)}\\
         &= \rho^{-2m_n}\sum_{k=1}^{N_n} \left(I_k(f_n)\right)^{2/(1+\rho^2)} + \rho^{-2m_n}\sum_{k>N_n} \left(I_k(f_n)\right)^{2/(1+\rho^2)}\\
         &\leq \rho^{-2m_n} {N_n}^{\rho^2/(1+\rho^2)} \left(\sum_{k=1}^{N_n} I_k(f)^2\right)^{1/(1+\rho^2)} + \rho^{-2m_n}\sum_{k>N_n} \left(I_k(f_n)\right)^{2/(1+\rho^2)}
    \end{split}
\end{equation*}

So, handling these terms separately, we see for the first that
$$\rho^{-2m_n} {N_n}^{\rho^2/(1+\rho^2)} \left(\sum_{k=1}^{N_n} I_k(f)^2\right)^{1/(1+\rho^2)} \leq \rho^{-2m_n} \left({N_n}^{\rho^2}H(f_n)\right)^{1/(1+\rho^2)}$$
and for the second, using \smartref{Corollary}{cor_influence_tail_bound} -- noting that we use the assumption that $p>\frac{1}{2}$ and are imposing a restriction on what $\rho$ can be chosen to be -- that
$$\rho^{-2m_n}\sum_{k>N_n} \left(I_k(f_n)\right)^{2/(1+\rho^2)} \leq \rho^{-2m_n}\mu_n^{2/(1+\rho^2)}N_n^{1-2p/(1+\rho^2)}\frac{1}{\frac{2p}{1+\rho^2}-1}$$
so that, putting them together again, we get that
$$\sum_{1 \leq \abs{S} \leq m_n} \hat{f_n}(S)^2 \leq \rho^{-2m_n}\left( \left({N_n}^{\rho^2}H(f_n)\right)^{1/(1+\rho^2)} + \frac{\mu_n^{2/(1+\rho^2)}N_n^{1-2p/(1+\rho^2)}}{\frac{2p}{1+\rho^2}-1}\right)$$

Now, if we apply our bound on $H(f_n)$ and take $N_n = \mu_n^\alpha$ for some $\alpha$ to be specified later, this becomes
\begin{equation*}
    \begin{split}
        \sum_{1 \leq \abs{S} \leq m_n} \widehat{f_n}(S)^2 &\leq \rho^{-2m_n}\left( \left({\mu_n}^{\alpha\rho^2}\mu_n^{-\gamma}\right)^{1/(1+\rho^2)} + \frac{\mu_n^{2/(1+\rho^2)}{\mu_n}^{\alpha(1-2p/(1+\rho^2))}}{\frac{2p}{1+\rho^2}-1}\right)\\
        &= \rho^{-2m_n}\left(\mu_n^{\frac{\alpha\rho^2 - \gamma}{1+\rho^2}} + \frac{\mu_n^{\alpha + (1-\alpha)\frac{2}{1+\rho^2}}}{\frac{2p}{1+\rho^2}-1}\right)
    \end{split}
\end{equation*}
and, taking $m_n$ to be logarithmic in $\mu_n$ and making appropriate choices of $\rho$ and $\alpha$, this can be made to go to zero as $n$ (and thus also $\mu_n$) goes to infinity. The details are left as an exercise for the reader's computer algebra system.
\end{proof}
\end{theorem}

%% file: 4_3_KKL.tex
\subsection{The Kahn-Kalai-Linial theorem on maximal influences}
In the case of a finitely-supported Boolean function, it makes sense to ask how small the largest influence for some bit can be, in relation to the total number of bits and the variance of the function. This of course does not transfer directly to the finitary setting, since we no longer have a fixed finite number of bits.

What does transfer easily from the finite case is the relationship between the total influence of a function and the maximal influence. The proof is basically exactly the same as in the finite case, applying hypercontractivity, except of course that we need hypercontractivity for finitary functions, not just finitely supported ones.

To begin with, however, we need a simpler result relating the total influence to the variance of a function:
\begin{lemma}[Poincaré's inequality]
Suppose $f$ is some Boolean function $\Omega^\infty\to\Omega$, not necessarily finitary. It then holds that
$$I(f) \geq \Var{f}$$
\begin{proof}
This is trivial from noticing that $\Var{f} = \sum_{\abs{S}>0} \widehat{f}(S)^2$ while \smartref{Lemma}{lemma_fourier_formula_for_influences} gives $I(f) = \sum_S \abs{S}\widehat{f}(S)^2$.
\end{proof}\end{lemma}

\begin{remark}In the finite setting, Poincaré's inequality can also be proved purely combinatorially, without resorting to Fourier formulas. In order to transfer this proof to the infinite-bit setting, we actually need to assume $f$ is finitary. We give the proof here, to give further illustration of how finitary functions behave like finitely-supported functions:
\begin{proof}
First off, note that since $f$ is Boolean we have, if $\omega$ and $\omega'$ are independent uniform, that
\begin{equation*}
    \begin{split}
        \Prob{f(\omega) \neq f(\omega')} &= \Prob{f(\omega)f(\omega') = -1}\\
        &= \frac{1 - \E{f(\omega)f(\omega')}}{2}\\
        &= \frac{\E{f^2} - \E{f}^2}{2} = \frac{\Var{f}}{2}
    \end{split}
\end{equation*}
so what we wish to do is to upper bound the probability that $f(\omega)\neq f(\omega')$.

So, for each $k\geq 0$, let $\omega_k = (\omega'_1, \omega'_2, \ldots, \omega'_{k-1}, \omega_k, \omega_{k+1},\ldots)$. That is, $\omega_k$ agrees with $\omega'$ up to bit $k-1$ and is then $\omega$ from bit $k$ onward.

Now, since $f$ is finitary, if $f(\omega) \neq f(\omega')$ there almost surely exists some $N$ such that $f(\omega) = f(\omega_0) \neq f(\omega_N) = f(\omega')$. We can also see that for $f(\omega_0) \neq f(\omega_N)$ to hold, there must exist some $k$ such that $f(\omega_k)\neq f(\omega_{k+1})$, since we are flipping bits one at a time. So, we have that
$$\left\{f(\omega) \neq f(\omega')\right\} \subseteq \bigcup_k \left\{f(\omega_k)\neq f(\omega_{k+1})\right\}$$
and thus a union bound gives
$$\Var{f} = 2\Prob{f(\omega)\neq f(\omega')} \leq 2\sum_k \Prob{f(\omega_k)\neq f(\omega_{k+1})} = \sum_k I_k(f)$$
\end{proof}

Notice that we actually needed finitariness to get that $f(\omega)\neq f(\omega')$ implies $f(\omega) \neq f(\omega_N)$ for some $N$. For a non-finitary function we could have the phenomenon that $f(\omega)\neq f(\omega')$, but no finite amount of flips to make $\omega$ match $\omega'$ actually changes the value of $f$. A non-finitary function is thus somewhat impervious to combinatorial reasoning.
\end{remark}

\begin{theorem}\label{theorem_KKL_v1}
There exists a universal constant $c>0$ such that if $f$ is any Boolean function, not necessarily finitary, and $\delta = \sup_i I_i(f)$, then
$$I(f) \geq c\Var{f}\log\left(\frac{1}{\delta}\right)$$

The same also holds if one extends the definition of $I_k(f)$ from Boolean $f$ to $f\in L^2$ using the Fourier formula of \smartref{Lemma}{lemma_fourier_formula_for_influences}.
\begin{proof}
We can divide into two cases -- $\delta \leq \frac{1}{1000}$ and $\delta > \frac{1}{1000}$. We start with the first case.

Take an integer $M$ to be specified later. We can then calculate, reusing some tricks from our proof of the BKS theorem, that
\begin{equation*}
    \begin{split}
        \sum_{1 \leq \abs{S} \leq M} \widehat{f}(S)^2 &\leq \sum_{1 \leq \abs{S} \leq M} \widehat{f}(S)^2 \abs{S}\widehat{f}(S)^2\\
        &\leq 2^{2M} \sum_{\abs{S}\geq 1} (1/2)^{2\abs{S}}\abs{S}\widehat{f}(S)^2\\
        &= \frac{1}{4}2^{2M}\sum_k \norm{T_{1/2}(\nabla_k f)}_2^2
    \end{split}
\end{equation*}
and obviously if $f\in L^2$ then $\nabla_k f$ will be in $L^2$, so we can apply hypercontractivity (\smartref{Theorem}{theorem_hypercontractivity}) to see that
\begin{equation*}
    \begin{split}
        \sum_{1 \leq \abs{S} \leq M} \widehat{f}(S)^2 &\leq \frac{1}{4}2^{2M}\sum_k \norm{\nabla_k f}_ {5/4}^2\\
        &\leq 2^{2M}\sum_k I_k(f)^{8/5}\\
        &\leq 2^{2M}\delta^{3/5} \sum_k I_k(f) = 2^{2M}\delta^{3/5}I(f)
    \end{split}
\end{equation*}

We can now see that
\begin{equation*}
    \begin{split}
        \Var{f} = \sum_{\abs{S}\geq 1} \widehat{f}(S)^2 & \leq \sum_{1\leq \abs{S} \leq M} \widehat{f}(S)^2 + \frac{1}{M}\sum_{\abs{S}>M} \abs{S}\widehat{f}(S)^2\\
        &\leq \left(2^{2M}\delta^{3/5} + \frac{1}{M}\right)I(f)
    \end{split}
\end{equation*}
so, making the obvious choice of $M = \frac{1}{2}\left(\frac{3}{5}\log_2\left(\frac{1}{\delta}\right) - \log_2\left(\log_2\left(\frac{1}{\delta}\right)\right)\right)$, and using that $\delta \leq 1/1000$ to see that $M \geq \frac{1}{10}\log_2\left(\frac{1}{\delta}\right)$, one can verify that we get
$$\Var{f} \leq \left(\frac{1}{\log_2\left(\frac{1}{\delta}\right)} + \frac{10}{\log_2\left(\frac{1}{\delta}\right)}\right)I(f)$$
which gives
$$I(f) \geq \frac{1}{11\log(2)}\Var{f}\log\left(\frac{1}{\delta}\right)$$
and so we have the result for $\delta \leq \frac{1}{1000}$.

For the case of $\delta>1/1000$, the discrete Poincaré inequality says that $I(f) \geq \Var{f}$, which gives us that in this case the claim is true with $c=\frac{1}{\log(1000)}$. This finishes the proof.
\end{proof}
\end{theorem}

The above result was entirely straightforward to carry over from the finite case, since we could just apply the same tools that we had already set up in our setting.

To carry over the $\frac{\log n}{n}$ bound to the finitary case, we clearly need to make some more drastic adjustments, since we no longer have an $n$. In the finite case, one crucial step is the observation that the maximal influence is certainly not less than the average influence, that is, that $\max_k I_k(f) \geq \frac{1}{n}I(f)$. So, we begin by transferring this argument to the infinite case, picking up an error term in the process.

\begin{lemma}\label{lemma_bound_maxinf_influence_over_nu}
Suppose $f$ is $p$-knowable for some $p>1$. It then holds that
$$\max_k I_k(f) \geq \frac{I(f)}{\E{\max W}} - \frac{1}{p-1}\frac{\E{\left(\max W\right)^p}}{\E{\max W}^p}$$
\begin{proof}
We calculate, applying \smartref{Corollary}{cor_influence_tail_bound}, that
\begin{equation*}
    \begin{split}
        \max_k I_k(f) &\geq \max_{k \leq \E{\max W}} I_k(f)\\
        &\geq \frac{1}{\E{\max W}} \sum_{k=1}^{\floor{\E{\max W}}} I_k(f)\\
        &= \frac{1}{\E{\max W}}I(f) - \frac{1}{\E{\max W}}\sum_{k>\E{\max W}} I_k(f)\\
        &\geq \frac{I(f)}{\E{\max W}} - \frac{1}{\E{\max W}}\frac{\E{\left(\max W\right)^p}\E{\max W}^{1-p}}{p-1}\\
        &= \frac{I(f)}{\E{\max W}} - \frac{1}{p-1}\frac{\E{\left(\max W\right)^p}}{\E{\max W}^p}
    \end{split}
\end{equation*}
as desired.
\end{proof}
\end{lemma}

Now this, together with Poincaré's inequality, immediately gives us a version of the easy $\frac{1}{n}$ bound from the finite case:
\begin{cor}
Suppose $f$ is $p$-knowable for some $p>1$. It then holds that
$$\max_k I_k(f) \geq \Var{f}\frac{1}{\E{\max W}} - \frac{1}{p-1}\frac{\E{\left(\max W\right)^p}}{\E{\max W}^p}$$
\end{cor}

We can, however, actually do better than this, and make the same logarithmic gain as in the finite case. To do this, we leverage our existing extension of the KKL theorems.

\begin{theorem}
Suppose $f$ is $p$-knowable for some $p>1$. There exists a universal constant $c$, that does not depend on $p$, such that
$$\max_k I_k(f) \geq c\Var{f}\frac{\log\left(\E{\max W}\right)}{\E{\max W}} - \frac{1}{p-1}\frac{\E{\left(\max W\right)^p}}{\E{\max W}^p}$$
\begin{proof}
Let $\delta = \max_k I_k(f)$, and $\nu = \E{\max W}$. We have two cases -- either $\delta \geq \nu^{-1/2}$ or $\delta < \nu^{-1/2}$. Clearly, in the first case, $\frac{1}{\sqrt{\nu}}$ is asymptotically greater than $\frac{\log(\nu)}{\nu}$, and so some appropriate constant exists.

In the second case, we can use our existing KKL theorem, \smartref{Theorem}{theorem_KKL_v1}, together with \smartref{Lemma}{lemma_bound_maxinf_influence_over_nu} to get that
\begin{equation*}
    \begin{split}
        \max_k I_k(f) &\geq \frac{I(f)}{\nu} - \frac{1}{p-1}\frac{\E{\left(\max W\right)^p}}{\E{\max W}^p}\\
        &\geq \frac{c\Var{f}\log\left(\frac{1}{\delta}\right)}{\nu} - \frac{1}{p-1}\frac{\E{\left(\max W\right)^p}}{\E{\max W}^p}\\
        &\geq \frac{c\Var{f}\log\left(\frac{1}{\nu^{-1/2}}\right)}{\nu} - \frac{1}{p-1}\frac{\E{\left(\max W\right)^p}}{\E{\max W}^p}\\
        &= \frac{\frac{1}{2}c\Var{f}\log(\nu)}{\nu} - \frac{1}{p-1}\frac{\E{\left(\max W\right)^p}}{\E{\max W}^p}
    \end{split}
\end{equation*}
as desired.

So, taking the least of the two constants from the two cases, the theorem is proven.
\end{proof}
\end{theorem}

\begin{remark}
Note that there is nothing here that restricts us to taking $W$ to be the least witness set -- all the arguments go through as long as $W$ is almost surely a finite witness set. So in the case of a function $f:\Omega^n\to\Omega$, we can take $\max W \equiv n$, and recover the usual finite case of the KKL theorem by taking the limit as $p$ goes to infinity in the finitary version.

Thus, the theorem we have proven here is at least genuinely a generalisation of the result from the finite case. How useful it actually is will of course depend on the behaviour of the error term that has appeared.
\end{remark}

%% file: 5_revealment_theorem.tex
\section{Algorithms, and the revealment theorem}

So far, we've seen the definition of noise sensitivity, and results connecting it with the spectrum of the functions and with their influences. It remains to carry over the third main condition for noise sensitivity in the finite case, the one in terms of algorithms.

\begin{definition}
For a finitary function $f: \Omega^\infty \to \Omega$, a (randomised) \emph{algorithm} $A$ looks at the bits of $\omega$ in some order -- which may depend on the bits it has seen so far and on some auxiliary randomness -- until it knows the value of the function. We assume that the algorithm is also finitary, in that the set of bits it looks at is almost surely finite.

Let the random set of bits looked at by $A$ be $W(A)=W(A)(\omega,\Tilde{\omega})$, where $\Tilde{\omega}$ is the randomness of the algorithm. As the notation implies, $W(A)$ will of course always be a witness set for $f$.
\end{definition}

\begin{prop}\label{prop_finitary_equiv_exists_algorithm}
A function $f$ is finitary if and only if there exists an algorithm for it that almost surely looks at only finitely many bits.
\begin{proof}
The set of bits looked at by an algorithm is always a witness set, since by definition the algorithm knows the value of the function when it stops, and looking at the bits in order until we have seen a witness set will always look at only finitely many bits if there is a finite witness set.
\end{proof}
\end{prop}

\begin{definition}
Given a finitary function $f$ and an algorithm $A$ for $f$, we define the \emph{revealment of $A$ on bit $i$} to be $\delta_i^A = \Prob{i \in W(A)}$, and the \emph{revealment of $A$} to be $\sup_i \delta^A_i$. Finally, we define the \emph{revealment of $f$} to be
$$\delta_f = \inf_A \delta^A = \inf_A \sup_i \Prob{i \in W(A)}$$
where the infimum runs over all algorithms $A$ for $f$.
\end{definition}

In the finite setting, we have the following theorem:
\begin{theorem}[\cite{schramm_steif_2010}]
\label{theorem_finite_case_revealments_ineq}
For any function $f: \Omega^n \to \R$ and for each $k\in\N$, we have that
$$E_f(k) = \sum_{S\subseteq[n], \abs{S}=k} \hat{f}(S)^2 \leq \delta_f k\norm{f}_2^2$$
\end{theorem}

From \smartref{Theorem}{theorem_finite_case_revealments_ineq} one easily gets, using \smartref{Theorem}{theorem_fourierclassificationofnoisesensitivity}, a criterion for noise sensitivity:
\begin{cor}\label{cor_finite_case_revealments_to_zero_NS}
Let $f_n: \Omega^{N_n} \to \Omega$ be a sequence of finitely supported Boolean functions. If
$$\delta_{f_n} \to 0$$
then the sequence $f_n$ is noise sensitive.
\end{cor}

In order to extend these results to the finitary setting, we need to know that $\delta_f$ is well approximated by $\delta_g$ if $g$ is a finite approximation of $f$. Once this is established, our previous machinery will easily give us our desired results also in the finitary setting. 

\begin{theorem}\label{theorem_revealments_approximable}
Let $f$ be some  finitary Boolean function. It holds that $\delta_{B_m f} \to \delta_f$ as $m\to\infty$, where $B_m f$ is the Boolean approximation of $f$ defined in \smartref{Definition}{definition_boolean_approximation}.
\begin{proof}
For this proof, we proceed in two parts -- first we show that $\limsup_m \delta_{B_m f} \leq \delta_f$ by constructing an algorithm for $B_m f$ given one for $f$, and then we show that $\liminf_m \delta_{B_m f} \geq \delta_f$ by constructing an algorithm for $f$ given ones for each $B_m f$.

\emph{Part one:}

For each $n$, take an algorithm $A_n$ for $f$ such that $\delta^{A_n} < \delta_f+\frac{1}{n}$. Now, for each $m\in\N$, define an algorithm $A_{n,m}$ for $B_m f$ by:
\begin{enumerate}
    \item Run $A_n$ for $f$. If it terminates without ever attempting to look at a bit beyond $m$, then we know the value of $f$ from the first $m$ bits, and it thus equals the value of $B_m f$, which we now know.
    \item If $A_n$ requests to look at a bit beyond $m$, instead look at every bit between $1$ and $m$. This will determine $B_m f$, and so we can terminate.
\end{enumerate}

We can now calculate, for each $i\in [m]$, the revealment of $A_{n,m}$ on bit $i$ as follows
\begin{equation*}
    \begin{split}
        \delta_i^{A_{n,m}} &= \Prob{i \in W(A_{n,m})}\\
        &= \Prob{i \in W(A_{n,m}) \given W(A_n) \subseteq [m]}\Prob{W(A_n) \subseteq[m]}\\
        &\qquad + \Prob{i \in W(A_{n,m}) \given W(A_n) \not\subseteq [m]}\Prob{W(A_n) \not\subseteq[m]}
    \end{split}
\end{equation*}
and, conditional on $W(A_n) \subseteq[m]$, we have $W(A_{n,m}) = W(A_n)$, and so the first term is simply $\Prob{i \in W(A_n), W(A_n) \subseteq[m]}$. For the second term, we have that $W(A_{n,m}) = [m]$ conditional on $W(A_n)\not\subseteq[m]$, so it reduces to just $\Prob{W(A_n) \not\subseteq[m]}$. 

Putting this together, we have the following
\begin{equation*}
    \begin{split}
        \delta_i^{A_{n,m}} &= \Prob{i \in W(A_n), W(A_n) \subseteq[m]} + \Prob{W(A_n) \not\subseteq[m]}\\
        &\leq \delta_i^{A_n} + \Prob{\max W(A_n) > m}
    \end{split}
\end{equation*}
and thus
$$\delta_{B_m f} \leq \delta^{A_{n,m}} \leq \delta^{A_n} + \Prob{\max W(A_n) > m} < \delta_f + \frac{1}{n} + \Prob{\max W(A_n) > m}$$

From this, we see that $\limsup_m \delta_{B_m f} \leq \delta_f$. Specifically, for each $\epsilon>0$, take an $n$ large enough that $\frac{1}{n}<\frac{\epsilon}{2}$. For this fixed $n$, we know that $\max W(A_n)$ is almost surely finite, and so in particular $\Prob{\max W(A_n) > m}$ is decreasing towards $0$ as $m$ increases. 

So, for large enough $m$, we have that $\Prob{\max W(A_n) > m} < \frac{\epsilon}{2}$. Putting this together, we have that $\delta_{B_m f} < \delta_f + \epsilon$ for all large enough $m$, showing our claim.

\emph{Part two:}

For each $n$, take for every $m$ an algorithm $A_{n,m}$ for $B_m f$ such that $\delta^{A_{m,n}} < \delta_{B_m f} + \frac{1}{n}$. Given this, we can for each pair $n, m$ construct an algorithm $C_{n,m}$ for $f$ as follows:
\begin{enumerate}
    \item Run $A_{n,m}$ for $B_m f$, and consider the set of bits looked at by the algorithm. If this is a witness set for $f$ (which can be determined by looking only at the bits so far, so this step requires no additional queries), then we know that $B_m f = f$ on this input and thus we also know the value of $f$ on this input, and can terminate.
    \item If the set of bits looked at by $A_{n,m}$ is not a witness set for $f$, look at the bits in order from the first and forwards until the value of $f$ is determined.
\end{enumerate}

Again, we can write the revealment of $C_{n,m}$ on bit $i$ as
\begin{equation*}
    \begin{split}
        \delta_i^{C_{n,m}} &= \Prob{i \in W(C_{n,m})}\\
        &= \Prob{i \in W(C_{n,m}) \given W(C_{n,m}) \subseteq[m]}\Prob{W(C_{n,m}) \subseteq[m]}\\
        &\qquad + \Prob{i \in W(C_{n,m}) \given W(C_{n,m}) \not\subseteq[m]}\Prob{W(C_{n,m}) \not\subseteq[m]}
    \end{split}
\end{equation*}
Now, by construction, conditional on $W(C_{n,m}) \subseteq[m]$ we have $W(C_{n,m}) = W(A_{n,m})$, and so the first term is $\Prob{i \in W(A_{n,m}), W(C_{n,m}) \subseteq[m]}$. Given this, we thus have that
$$\delta_i^{C_{n,m}} \leq \Prob{i \in W(A_{n,m})} + \Prob{W(C_{n,m}) \not \subseteq[m]}$$

For the second term, we have by construction that $W(C_{n,m}) \not\subseteq[m]$ precisely when there's \emph{no} subset of $[m]$ that is a witness set. Thus, letting $W(f)$ be our minimal witness set as in the very definition of being finitary, \smartref{Definition}{definition_finitary_function}, the second term is precisely $\Prob{\max W(f) > m}$. So, taking a supremum over all $i$, we arrive at that
$$\delta^{C_{n,m}} \leq \delta^{A_{n,m}} + \Prob{\max W(f) > m}$$

Now, since $C_{n,m}$ is an algorithm for $f$, we have $\delta_f \leq \delta^{C_{n,m}}$, and by our choice of $A_{n,m}$ we have $\delta^{A_{n,m}} < \delta_{B_m f} + \frac{1}{n}$. So our previous inequality becomes
$$\delta_f < \delta_{B_m f} + \frac{1}{n} + \Prob{\max W(f) > m}$$
which, letting $n\to\infty$ and observing that $\Prob{\max W(f) > m}$ decreases towards $0$ with $m$ since $f$ is finitary, shows our claim that $\liminf \delta_{B_m f} \geq \delta_f$.

So we have shown that
$$\delta_f \leq \liminf \delta_{B_m f} \leq \limsup \delta_{B_m f} \leq \delta_f$$
so the limit exists and equals $\delta_f$, as desired.
\end{proof}
\end{theorem}

With this result in hand, giving our extensions to the finitary case is just an exercise in cross-referencing.

\begin{cor}\label{cor_finitary_revealment_ineq}
For any finitary Boolean function $f$ and for each $k\in\N$, we have that
$$E_f(k) = \sum_{S\subset\N, \abs{S} = k} \hat{f}(S)^2 \leq \delta_f k \norm{f}_ 2^2$$
\begin{proof}
Apply \smartref{Theorem}{theorem_finite_case_revealments_ineq} to $B_n f$, and then take the limit as $n\to\infty$. \smartref{Theorem}{theorem_revealments_approximable} gives $\delta_{B_n f} \to \delta_f$, and \smartref{Proposition}{prop_boolean_approximations_converge} gives that $B_n f \to f$ in $L^2$. Together, this gives that we converge to the right inequality.
\end{proof}
\end{cor}

\begin{cor}\label{cor_finitary_revealments_to_zero_NS}
If $f_n$ is a sequence of finitary Boolean functions such that $\delta_{f_n} \to 0$, then $f_n$ is noise sensitive.
\begin{proof}
This can either be seen through \smartref{Corollary}{cor_finitary_revealment_ineq} and \smartref{Theorem}{theorem_fourierclassificationofnoisesensitivity}, or by using \smartref{Theorem}{theorem_revealments_approximable} directly together with \smartref{Corollary}{cor_finite_case_revealments_to_zero_NS} and \smartref{Proposition}{prop_approximations_noise_sensitive_implies_noise_sensitive}.
\end{proof}
\end{cor}

%% file: 6_example_voter_model.tex
\section{Example: Discrete-time edge-ordered voter model}

Having set up a lot of theory, we finally arrive at an application to show how we can use the theory we have established. The setting is essentially the voter model, except that we have to adjust it a bit to be able to encode it using a Boolean function. For a treatment of the usual voter model, including the duality we exploit here, see for example \cite{grimmett_2018}. Since this is just an example, we will handwave many of the details of the model -- all the properties we claim for our model are easily seen by the same method as for the usual voter model.

Consider some directed graph $G = (V,E)$, with an ordering of the edges $e_1, e_2, \ldots, e_n$. We take some subset $A_0 \subseteq V$ of the vertices which we consider to be coloured black at time zero, the rest being white. The discrete-time synchronous voter model evolves by that, at each time step, the next edge in order transfers colour along it with probability $0.5$, and otherwise does nothing. Once we reach time $n+1$, we loop around the list and edge $e_1$ is considered, and so on.

More formally, given an edge $e_i$, denote its origin by $e_i^o \in V$ and its target by $e_i^t \in V$. Given a state $A_{k-1} \subseteq V$ at time $k-1$, we generate $A_{k}$ as follows:
\begin{enumerate}
    \item With probability $\frac{1}{2}$, $A_k = A_{k-1}$.
    \item Otherwise, let $j = k \mod n$. If $e_j^o \in A_{k-1}$, set $A_k = A_{k-1}\cup\{e_j^t\}$, and if $e_j^o \not\in A_{k-1}$, set $A_k = A_{k-1} \setminus \{e_j^t\}$.
\end{enumerate}

Now, if we assume that $G$ is strongly connected, it is easily seen that this model will almost surely fixate in either the state $A_k = V$ or $A_k = \emptyset$ at some finite time $k$. The problem we are interested in is which of these two will occur.

It is easy to see that we can, for each specific $v\in V$, determine its colour at time $k$ by simply tracking backwards along the edges it got its colour from until we reach time $0$. This then becomes a random walk going in the opposite time-direction. Formally, we need to define this backward walk as taking one step for each full cycle through the edges, in order for it to be a time-homogeneous Markov chain. If we track two of these walks for two vertices $v, v' \in V$ backwards and they meet at some point, then they will clearly stay together from that point onward.

Note also that, since we assume the graph is finite and strongly connected, and there is a positive probability of staying put, the backwards walk will have a stationary distribution. Call this distribution $\pi_G$, remembering that the ordering of the edges is part of the data in $G$.

\begin{figure}
  \caption{Example trajectory of edge-directed voter model on a line graph with five vertices. Orange edges are activated, black ones are not. Time evolves from left to right. Colouring at time zero is omitted, since it is essentially irrelevant -- what we really want to study is which vertex the fixating colour comes from.\label{fig_example_voter_model}}
  \centering
    \includegraphics[width=0.7\textwidth]{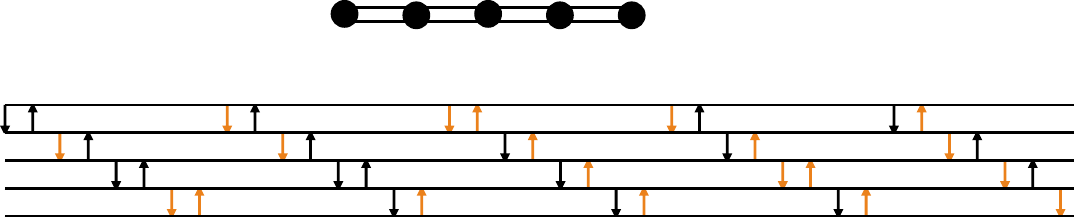}
\end{figure}

So, if we wish to determine which state the model fixates in, we can proceed as follows: Pick some very large time $N$, at which it is very likely to have fixated. For each vertex, track backwards -- if our family of random walks coalesces, then we know that the model has fixated, and can continue tracking backwards to time zero to find which state it fixated in. If this fails, we can just try again with an even larger $N$.

This method for determining the state it fixated in suggests that the initial state is not very relevant -- given just the randomness and not the initial state, we can determine a vertex $\rho\in V$ such that (given the randomness) which state the model fixates in is entirely determined by whether $\rho \in A_0$ or $\rho \in A_0^c$.

Here, we see an opening to insert our theory of Boolean functions -- encoding the choices of randomness for each edge at each time as an element of $\Omega^\infty$, this $\rho$ becomes a function from $\Omega^\infty$ to $V$. That this will be a finitary function follows from that the model almost surely fixates.

\begin{figure}
  \caption{Backwards walk trajectory associated to the example trajectory of \smartref{Figure}{fig_example_voter_model}. The orange lines give the trajectory of the backwards-tracking walk, which as can be seen coalesces and then does a random walk until time zero. Note that, reading left to right, we get an edge-oriented voter model trajectory with each edge corresponding to one time-step, while reading right-to-left we get a backwards walk trajectory with \emph{each group of eight edges} corresponding to one time-step.}
  \centering
    \includegraphics[width=0.7\textwidth]{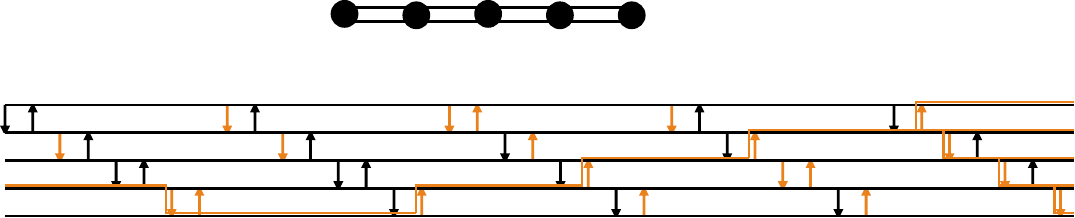}
\end{figure}

Now, if we refine our description of the algorithm for $\rho$ above, we can get noise sensitivity for the voter model. Effectively, the problem with the algorithm above is that we will at some point in time be asking for every single bit, to determine the first few steps of the coalescing walk process. Once the walk has coalesced, and the single walker has become very nearly distributed according to the stationary distribution, we have a good quantity to bound by. The issue is getting rid of the times before then.

The way we solve this problem is like with so many other problems -- we procrastinate and put it off to some vague far off future. It will turn out that if we just let it be far off enough, and are vague enough about when exactly we will do it, the problem goes away on its own.

\begin{theorem}
Let $G_n = (V_n, E_n)$ be some family of strongly connected graphs with an ordering of the edges, and let $B_n \subseteq V_n$ be some family of subsets of the vertices of each graph. For each graph, define $\rho_n: \Omega^\infty \to V_n$ as in the discussion above.

Let, for each $G_n$, for all $v, w \in V_n$, $\pi_k(v,w)$ be the amount of paths from $v$ to $w$ such that the sequence of edges is increasing with respect to the numbering of the edges. In order to get noise sensitivity, we will need that this quantity does not grow exponentially in $k$. In particular, define for every graph $G$ with ordered edges
$$\zeta(G) = \sup_{v\in V}\sum_{w\in V}\left(\sum_{k=1}^\infty\frac{\abs{\pi_k(v,w)}}{2^k}\right)$$

If the $G_n$ are spread out in the sense that
$$\lim_{n\to\infty} \left(\max_{v\in V_n} \pi_{G_n}(v)\right)\zeta(G_n) = 0$$
then the sequence of functions
$$2\ind{\rho_n \in B_n}-1: \Omega^\infty \to \Omega$$
is noise sensitive.
\begin{proof}
Note that $\zeta(G_n)$ is always at least $\frac{1}{2}$, so our hypothesis implies in particular that $\max_{v\in V} \pi_{G_n}(v)$ goes to zero. Thus, it is clear from \smartref{Corollary}{cor_finitary_revealments_to_zero_NS} that if we can, given a strongly connected graph $G=(V,E)$ with an ordering of the edges, show that $\delta_\rho \leq \left(\max_{v\in V} \pi_{G}(v)\right)(1+2\zeta(G))$, the result will follow.

So, take some such graph $G$, and define given integers $N$ and $M$ an algorithm $A_{N,M}$ for $\rho$ as follows:
\begin{enumerate}
    \item Using the algorithms auxiliary randomness, choose a random integer $U$ uniformly on ${0,1,\ldots, M}$.
    \item At time $N+U$, start tracking backwards from every vertex. So, at each time step, for each random walk, if the random walk is at $v\in V$, query each edge $e: v' \to v$ whether it gave colour to $v$ at this time step. If it did, ask every edge $e': v'' \to v'$ preceding $e$ whether it gave colour to $v'$. Continue this process for each $v''$ you got a ``yes'' from until you either get all ``no''s or you are out of edges to query. This determines what step the backward random walks take.
    \item If the random walks coalesce, then we know that the model has fixated at time $N+U$, and we find $\rho$ by tracking back to time zero.
    \item If they fail to coalesce, just look at bits in order from zero and onward until $\rho$ has been determined.
\end{enumerate}
    
We now wish to compute the revealments of this algorithm -- so, fix some time $i\in\N$ for the backward walk\footnote{Recall that the time for the voter model is one time-step per edge, while the backward walk takes one time step per entire cycle of the walk.} and some edge $e_j \in E$ -- we wish to compute the revealment of the bit $b = in + j$ coding for edge $e_j$ at time $i$. 

Let $S_{N,M}$ be the event that the random walks going backwards from $N+U$ coalesce before reaching time $N$.\footnote{Recall that time goes backwards here -- this means it coalesces at some time between $N$ and $N+U$.} So we have
\begin{equation*}
    \begin{split}
        \delta_b^{A_{N,M}} &= \Prob{b\in W(A_{N,M})}\\
        &\leq \Prob{b \in W(A_{N,M}), S_{N,M}} + \Prob{S_{N,M}^c}
    \end{split}
\end{equation*}
and the fact that the random walks almost surely coalesce in finite time gives us that $\Prob{S_{N,M}^c}=o(1)$ as $M\to\infty$.

Now let $C$ be the time it takes for the random walks to coalesce, so that they coalesce at time $N+U-C$, and pick some $\epsilon>0$. Since the backward random walk has a finite state space and its distribution converges to the stationary distribution, we can take some $K_\epsilon\in\N$ such that, for any $k\geq K_\epsilon$, the distribution of the backward random walk at time $k$ is within $\epsilon$ of $\pi_G$ in the supremum norm, independently of the starting distribution.

So, continuing our calculation, we have
\begin{equation*}
    \begin{split}
        \delta_b^{A_{N,M}} &\leq \Prob{b \in W(A_{N,M}), S_{N,M}} + o(1)\\
        &= \Prob{b \in W(A_{N,M}), S_{N,M}, i \in \{1,\ldots, N+U-C-K_\epsilon\}}\\
        &\qquad + \Prob{b \in W(A_{N,M}), S_{N,M}, i \in \{N+U-C-K_\epsilon+1, \ldots, N+U\}}\\
        &\qquad + \Prob{b \in W(A_{N,M}), S_{N,M}, i > N+U} + o(1)
    \end{split}
\end{equation*}
where the final probability in the sum is obviously exactly zero, since the algorithm looks at no bits at times after $N+U$ if the random walks do coalesce. We claim that the second is also $o(1)$.

So, pick some $\eta>0$. Since $G$ is fixed and the distribution of $C$ depends only on $G$, we can pick some fixed integer $L$ such that $\Prob{C>L}<\frac{\eta}{2}$. So, denoting the second term in the sum by $*$, we can calculate
\begin{equation*}
    \begin{split}
        * &\leq \Prob{ S_{N,M}, i \in \{N+U-C-K_\epsilon+1, \ldots, N+U\}}\\
        &= \Prob{ S_{N,M}, i \in \{N+U-C-K_\epsilon+1, \ldots, N+U\}, C \leq L}\\
        &\qquad + \Prob{ S_{N,M}, i \in \{N+U-C-K_\epsilon+1, \ldots, N+U\}, C > L}\\
        &< \Prob{ S_{N,M}, i \in \{N+U-C-K_\epsilon+1, \ldots, N+U\}, C \leq L} + \frac{\eta}{2}\\
        &\leq \Prob{i \in \{N+U-L-K_\epsilon+1, \ldots, N+U\} \given S_{N,M}, C \leq L} + \frac{\eta}{2}
    \end{split}
\end{equation*}
Now, this probability is precisely the probability of some fixed integer lying in a random interval of length $L+K_\epsilon$ in $[M]$ whose endpoint is uniformly distributed. This is easily seen to go to zero with $M$, and so in particular can be gotten smaller than $\frac{\eta}{2}$. So what we have shown is that $*<\eta$ for large enough choices of $M$, which establishes our claim.

In summary, what we have seen so far is that
\begin{equation*}
    \begin{split}
        \delta_b^{A_{N,M}} &\leq \Prob{b \in W(A_{N,M}), S_{N,M}, i \in \{1,\ldots, N+U-C-K_\epsilon\}} + o(1)\\
        &\leq \Prob{b \in W(A_{N,M}) \given S_{N,M}, i \in \{1,\ldots, N+U-C-K_\epsilon\}} + o(1)
    \end{split}
\end{equation*}
and we are finally at the point of analysing the actual algorithm. So, from now on we assume we have always conditioned on the backwards walks having coalesced for a long enough time that the distribution of the walk is within $\epsilon$ of being stationary.

So, how can bit $b$ come to be queried? It can either code for an edge pointing at the current location of the random walk, or it can code for an edge pointing at a vertex which in turn gave its colour to the current location, or it can be three steps of colour-transfer away, and so on. Additionally, the paths connecting it to the current location need to have increasing numbers.

To clarify the problem we are attacking, denote the target of the edge which $b$ codes for by $v$, and the random position of the backward walk at the time by $X$. What we are interested in is the probability that there exists a path from $v$ to $X$, such that all except possibly the first edge in it transferred colour in the current cycle, and such that the edges form an increasing sequence with respect to the numbering of the edges.

Denote the collection of all such increasing paths from $v$ to $X$ by $\pi(v,X)$, including the trivial empty path in the case of $v=X$. Denote the subcollection of paths of length $k$ by $\pi_k(v,X)$. What we are interested in is the probability that there exists a path $\mathfrak{p} \in \pi(v,X)$ such that all except possibly the first edge in it are active. Denote this property of a path by $\mathcal{O}(\mathfrak{p})$.

With all of this notation set up, we can finally calculate that
\begin{equation*}
    \begin{split}
        \Prob{\exists \mathfrak{p} \in \pi(v,X): \mathcal{O}(\mathfrak{p})} &= \sum_{w \in V}\Prob{\exists \mathfrak{p} \in \pi(v,w): \mathcal{O}(\mathfrak{p}) \given X = w}\Prob{X=w}\\
        &= \sum_{w \in V}\Prob{\bigcup_{k=0}^\infty \bigcup_{\mathfrak{p} \in \pi_k(v,w)} \mathcal{O}(\mathfrak{p})}\Prob{X=w}\\
        &\leq \sum_{w\in V}\left(\sum_{k=0}^\infty\sum_{\mathfrak{p} \in \pi_k(v,w)}\Prob{\mathcal{O}(\mathfrak{p})}\right)\Prob{X=w}\\
        &= \sum_{w\in V}\left(\abs{\pi_0(v,w)} + \sum_{k=1}^\infty\sum_{\mathfrak{p} \in \pi_k(v,w)}2^{-(k-1)}\right)\Prob{X=w}\\
        &= \Prob{X=v} + \sum_{w\in V}\left(\sum_{k=1}^\infty\frac{\abs{\pi_k(v,w)}}{2^{k-1}}\right)\Prob{X=w}
    \end{split}
\end{equation*}
and now recall that the distribution of $X$ is within $\epsilon$ of the stationary distribution, so we can further compute
\begin{equation*}
    \begin{split}
        \Prob{\exists \mathfrak{p} \in \pi(v,X): \mathcal{O}(\mathfrak{p})} &\leq \Prob{X=v} + \sum_{w\in V}\left(\sum_{k=1}^\infty\frac{\abs{\pi_k(v,w)}}{2^{k-1}}\right)\Prob{X=w}\\
        &\leq \pi_G(v) + \epsilon + \sum_{w\in V}\left(\sum_{k=1}^\infty\frac{\abs{\pi_k(v,w)}}{2^{k-1}}\right)(\pi_G(w)+\epsilon)\\
        &\leq \left(\sup_{v'\in V} \pi_G(v') + \epsilon\right)\left(1 + \sum_{w\in V}\left(\sum_{k=1}^\infty\frac{\abs{\pi_k(v,w)}}{2^{k-1}}\right)\right)
    \end{split}
\end{equation*}

So, recalling where we started, what we have shown is that for every $\epsilon > 0$
$$\delta_b^{A_{N,M}} \leq \left(\sup_{v'\in V} \pi_G(v') + \epsilon\right)\left(1 + \sum_{w\in V}\left(\sum_{k=1}^\infty\frac{\abs{\pi_k(v,w)}}{2^{k-1}}\right)\right) + o_{N,M}(1)$$
and taking a supremum over $b$ in the left hand side means taking a supremum over $v$. So, taking that supremum, and then letting $\epsilon \to 0$ and $N, M \to \infty$, we finally get that
\begin{equation*}
    \begin{split}
        \delta_\rho &\leq \left(\sup_{v'\in V} \pi_G(v')\right)\left(1 + \sup_{v\in V}\sum_{w\in V}\left(\sum_{k=1}^\infty\frac{\abs{\pi_k(v,w)}}{2^{k-1}}\right)\right)\\
        &= \left(\sup_{v'\in V} \pi_G(v')\right)\left(1 + 2\zeta(G)\right)
    \end{split}
\end{equation*}
which is precisely the bound we were working towards.
\end{proof}
\end{theorem}

Now, this result establishes noise sensitivity for the voter model when the noise is applied to the randomness of the process. There are many other ways we could choose to phrase questions of noise sensitivity -- some of them have easy answers, while some seem very difficult. In addition to noising different things, we could also have different variants of the model -- our choice above of an edge-ordered discrete-time model was of course only to make our tools apply to it, not any preference for that model in particular.

Sticking with seeing the outcome of the process as being the state it fixates in, we could also apply the noise to the initial state. However, here, our discussion about the dual picture of a coalescing random walk easily translates into a proof that this will be noise stable. In fact, as is known, one can write
$$\Prob{\text{Fixate in black}} = \pi_G(A_0) = \sum_{v\in V} \pi_G(v)\ind{v\in A_0}$$
so the entire energy spectrum stays at level one.

If we were instead to apply the noise to the graph itself, there does not appear to be any easy answer to the question of noise sensitivity. Similarly, if we change the outcome from which state it fixates in to the time it takes to fixate, the problem appears to be difficult.